\title{A perturbative approach to the reconstruction of the eigenvalue
  spectrum of a normal covariance matrix from a spherically truncated
  counterpart\\[2.0ex]}   

\author{
  Filippo Palombi$^{\it a,b}$\footnote{Corresponding author. e--mail: 
  {\tt filippo.palombi$@$enea.it}}\ \ %
  and \  Simona Toti$^{\it a}$ \\[2.0ex]
  {$^{\it a}$\it \small ISTAT -- Istituto Nazionale di Statistica}\\[0.5ex]
  {\small Via Cesare Balbo 16, 00184 Rome -- Italy}\\[2.0ex]
  {$^{\it b}$\it \small ENEA -- Italian Agency for New Technologies, Energy}\\[0.0ex]
  {\it \small and Sustainable Economic Development}\\[0.5ex]
  {\small Via Enrico Fermi 45, 00044 Frascati -- Italy}\\[1.0ex]
}
\date{April 2014}

\documentclass[11pt,fleqn]{article}

\usepackage{footmisc,booktabs}
\usepackage{psfrag}
\usepackage{amsmath}
\usepackage{mathtools}
\usepackage{amsfonts}
\usepackage{dsfont}
\usepackage{amsthm}
\usepackage{epsfig}
\usepackage{multicol}
\usepackage[title]{appendix}
\usepackage{pgf}
\usepackage{bbm}
\usepackage[utf8]{inputenc}
\usepackage[margin=30pt,font=small,labelfont=bf,labelsep=endash]{caption}
\usepackage{marvosym}
\usepackage{lscape}
\usepackage{verbatim}
\usepackage{stackrel}
\usepackage{cite}
\usepackage{enumitem}

\numberwithin{equation}{section}\usepackage{amssymb}
\definecolor{Blue}{rgb}{0,0,1}
\definecolor{Black}{rgb}{0,0,0}

\newcommand{\bmu}{{\tilde\mu}}
\newcommand{\bk}{K}
\newcommand{\blambda}{{\tilde\lambda}}
\newcommand{\lambdaT}{{\lambda_{\rm\scriptscriptstyle T}}}
\newcommand{\lambdaEX}{{\lambda_{\rm ex}}}
\newcommand{\bR}{\tilde R}

\newcommand{\rhs}{{\it r.h.s.}\ }

\newcommand{\Oeone}{{\rm O}(\epsilon^1)}
\newcommand{\Oetwo}{{\rm O}(\epsilon^2)}
\newcommand{\Oethree}{{\rm O}(\epsilon^3)}
\newcommand{\Oefour}{{\rm O}(\epsilon^4)}
\newcommand{\Oen}{{\rm O}(\epsilon^n)}
\newcommand{\Oenpo}{{\rm O}(\epsilon^{n+1})}
\newcommand{\rL}{{\scriptscriptstyle\rm L}}
\newcommand{\rPT}{{\scriptscriptstyle\rm P}}
\newcommand{\Maple}{{\sc Maple}{\small\texttrademark}\ }

\newcommand{\loko}{\lambda_k^{(1)}}
\newcommand{\ltko}{\lambda_k^{(2)}}

\newcommand{\cB}{{\cal B}}

\newcommand{\cD}{{\cal D}}
\newcommand{\cE}{{\cal E}}

\newcommand{\cG}{{\cal G}}
\newcommand{\cK}{{\cal K}}
\newcommand{\cI}{{\cal I}}
\newcommand{\cJ}{{\cal J}}
\newcommand{\cM}{{\cal M}}
\newcommand{\cN}{{\cal N}}
\newcommand{\cO}{{\cal O}}
\newcommand{\cP}{{\cal P}}

\newcommand{\cT}{{\cal T}}

\newcommand{\rd}{{\rm d}}

\newcommand{\cov}{{\rm cov}}
\newcommand{\var}{{\rm var}}

\newcommand{\re}{{\rm e}}
\newcommand{\I}{{\mathbb{I}}}
\newcommand{\E}{{\mathbb{E}}}

\newcommand{\RR}{{\mathbb{R}}}

\newcommand{\trans}[1]{{#1}^{\scriptscriptstyle{\rm T}}}

\newcommand{\sign}{{\rm sign\,}}
\newcommand{\diag}{{\rm diag}}

\newcommand{\dR}{\mathds{R}}

\newcommand{\widesim}[2][3]{
  \mathrel{\overset{#2}{\scalebox{#1}[1.0]{$\sim$}}}
}

\newcommand{\cdf}{{\it c.d.f.}\ }
\newtheorem{prop}{Proposition}[section]
\newtheorem{corol}{Corollary}[section]

\makeatletter
\newcommand{\raisemath}[1]{\mathpalette{\raisem@th{#1}}}
\newcommand{\raisem@th}[3]{\raisebox{#1}{$#2#3$}}
\makeatother

\def\lambdabar{\protect\@lambdabar} 
\def\@lambdabar{%
\relax 
\bgroup 
\def\@tempa{\hbox{\raise.73\ht0 
\hbox to0pt{\kern.25\wd0\vrule width.5\wd0 
height.1pt depth.1pt\hss}\box0}}%
\mathchoice{\setbox0\hbox{$\displaystyle\lambda$}\@tempa}%
{\setbox0\hbox{$\textstyle\lambda$}\@tempa}%
{\setbox0\hbox{$\scriptstyle\lambda$}\@tempa}%
{\setbox0\hbox{$\scriptscriptstyle\lambda$}\@tempa}%
\egroup 
}

\definecolor{shadecolor}{rgb}{0.97,0.97,0.97}

\begin{document}
\maketitle

\begin{abstract}
In this paper we propose a perturbative method for the reconstruction of the covariance matrix of a multinormal distribution, under the assumption that the only available information amounts to the covariance matrix of a spherically truncated counterpart of the same distribution. We expand the relevant equations up to the fourth perturbative order and discuss the analytic properties of the first few perturbative terms. We finally compare the proposed approach with an exact iterative algorithm (discussed in Palombi et al. (2017)) in the hypothesis that the spherically truncated covariance matrix is estimated from samples of various sizes. 
\end{abstract}

\section{Introduction}

The analysis of the probability content of the multivariate normal distribution in finite (or infinite) regions has become a standard problem in probability and statistics since the original works by Ruben, refs.~\cite{ruben,ruben2,ruben3,ruben4}. Other aspects of the problem have been subsequently considered in refs.~\cite{tallis,tallis2} and more recently in ref.\cite{horrace}. In a recent paper \cite{palombi4} we studied how the covariance matrix $({\frak S}_\cB)_{ij}=\cov\left(X_i,X_j\,|\,X\in\cB_v(\rho)\right)$ of a multinormal random vector $X = \{X_k\}_{k=1}^v\sim\cN_v(0,\Sigma)$ in $v\ge 1$ dimensions, conditioned to a centered spherical domain $\cB_v(\rho)=\{x\in\RR^v:\ \trans{x}x \le \rho\}$, relates to the unconditioned covariance matrix $\Sigma_{ij}=\cov(X_i,X_j)$. In the same paper we also described practical situations where this kind of distributional 
truncation occurs.

\begin{center}
  \begin{figure}[!t]
    \centering
    \includegraphics[width=0.36\textwidth]{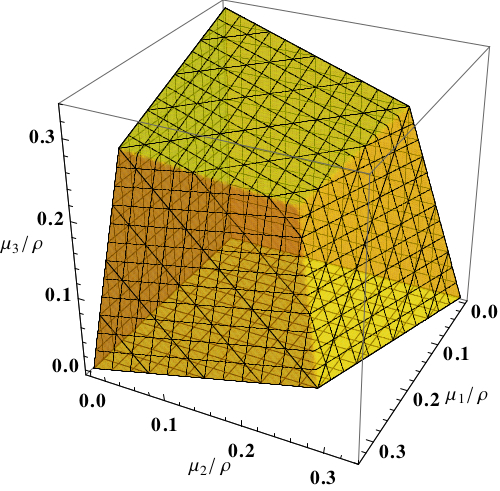}
    \hskip 2.0cm
    \includegraphics[width=0.36\textwidth]{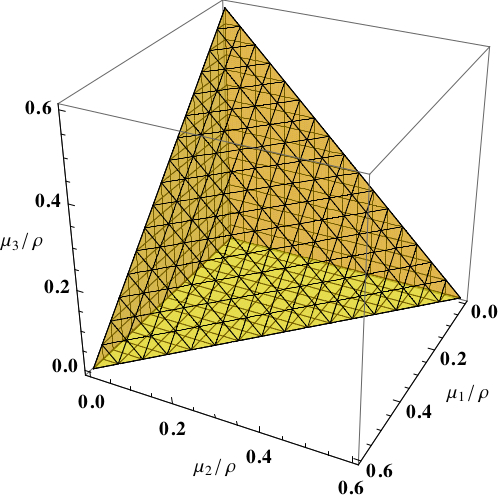}
    \vskip 0.2cm
    \caption{\small Domains $\cD(\tau_\rho^{-1})$ (left) and $\cD(\cT_\rho^{-1})$ (right) for $v=3$.\label{fig:defdoms}}    
    \vskip -0.3cm
  \end{figure}
\end{center}
\vskip -1.0cm
Owing to the symmetries of the geometrical set-up, ${\frak S}_\cB$ and $\Sigma$ can be shown to commute, i.e. if we let $\mu=\{\mu_k\}_{k=1}^v$ and $\lambda = \{\lambda_k\}_{k=1}^v$ denote respectively the eigenvalues of ${\frak S}_\cB$ and $\Sigma$, then  ${\frak S}_\cB = R\,\diag(\mu)\,\trans{R}$ and $\Sigma = R\,\diag(\lambda)\,\trans{R}$, with $R$ the same diagonalizing matrix. The eigenvalues fulfill the equations
\begin{equation} 
\mu_k =
\lambda_k\frac{\alpha_k(\rho;\lambda)}{\alpha(\rho;\lambda)}\,,\qquad\qquad k = 1,\ldots,v\,,
\label{eq:truncspectrum}
\end{equation}
with $\alpha$ and $\alpha_k$ belonging to the class of Gaussian integrals
\begin{equation}
\alpha_{k\ell m\dots}(\rho;\lambda) = \int_{\cB_v(\rho)}\rd^v x\
\frac{x_k^2}{\lambda_k}\,\frac{x_\ell^2}{\lambda_\ell}\,\frac{x_m^2}{\lambda_m}\dots\,\prod_{j=1}^v\delta(x_j,\lambda_j)\,,
\qquad \delta(y,\eta)\, =\,
\frac{\,\re^{-y^2/({2\eta})}\,}{(2\pi\eta)^{1/2}}\,.
\label{eq:alphaints}
\end{equation}
Since ${\frak S}_\cB$ and $\Sigma$ are simultaneously diagonalizable, we can assume  ${\frak S}_\cB = \diag(\mu)$ and  $\Sigma=\diag(\lambda)$ with no loss of generality. We are interested in reconstructing $\Sigma$ from ${\frak S}_\cB$, namely we aim at solving eqs.~(\ref{eq:truncspectrum}) with respect to $\lambda$ under the assumption that $\rho$ and $\mu$ are given (if ${\frak S}_\cB$ and $\Sigma$ are not diagonal and we assume that ${\frak S}_\cB$ is known, we can obtain $R$ from ${\frak S}_\cB$ by diagonalization and use it then for $\Sigma$.).  As explained in ref.~\cite{palombi4}, such reconstruction can be of practical importance, for example, in compositional analysis of multivariate log-normal data affected by outlying contaminations. Reconstructing $\Sigma$ from ${\frak S}_\cB$ is equivalent to writing $\lambda = \tau_\rho^{-1}\cdot \mu$, where $\tau_\rho^{-1}$ is a non-linear operator performing the covariance reconstruction. The existence of $\tau_\rho^{-1}$ follows from the injectivity of $\tau_\rho$. For the sake of completeness, we prove this in App.~A. The operator $\tau_\rho^{-1}$ is defined within the convex domain 
\begin{align}
  \cD(\tau_\rho^{-1}) \,\equiv \text{Im}(\tau_\rho) = \ & \bigcap_{k=1}^v\biggl\{\mu\in\dR_+^v: \ \sum_{j\ne k}^{1,\ldots,v} \mu_j + 3\mu_k<\rho\biggr\}\,.
\label{eq:domtauromo}
\end{align}
Eq.~(\ref{eq:domtauromo}) was proposed in ref.~\cite{palombi4} as a result of extensive numerical tests. In App.~B, we give analytic arguments to support its correctness, although we lack a complete proof. In Fig.~\ref{fig:defdoms}~(left), we plot $\cD(\tau_\rho^{-1})$ for $v=3$ as an illustrative example.

Unfortunately, there exists no general closed-form solution to eqs.~(\ref{eq:truncspectrum}), due to the non-linear character of the problem. For this reason, in ref.~\cite{palombi4} we proposed a numerical technique based on a fixed point iteration, whose convergence is ensured by some non-trivial variance and covariance inequalities within $\cB_v(\rho)$ \cite{palovar,mukerjee}.   

In this paper, we approach eqs.~(\ref{eq:truncspectrum}) from a different standpoint. We move from the observation that a simplified set-up occurs when the eigenvalues are fully degenerate, a case that was first considered by Tallis in ref.~\cite{tallis}. If $\mu_1 = \ldots = \mu_v \equiv \bmu$, by symmetry it follows $\lambda_1=\ldots = \lambda_v\equiv \blambda$ and the other way round. Eqs.~(\ref{eq:truncspectrum}) reduce in this limit to  
\begin{equation}
\bmu = \blambda\frac{F_{v+2}}{F_v}\biggl(\frac{\rho}{\blambda}\biggr)\, \equiv \,\cT_\rho(\blambda)\,,
\label{eq:tspectrum}
\end{equation}
with $F_v(x)$ denoting the \cdf of a $\chi^2$-variable with $v$ degrees of freedom\footnote{With abuse of notation we shall write $\frac{F_{v+2}}{F_v}(\frac{\rho}{\blambda})$ in place of $F_{v+2}(\rho/\blambda)/F_{v}(\rho/\blambda)$.}. It can be easily checked that $\cT_\rho(\blambda)$ is a monotonic increasing function of $\blambda$. In addition, $\cT_\rho(\blambda)$ fulfills
\begin{equation}
i{\text )}\ \,\lim_{\blambda\to 0}\cT_\rho(\blambda) = 0\,,\qquad\qquad
ii{\text )}\ \,{\lim_{\blambda\to \infty}\cT_\rho(\blambda) = \frac{\rho}{v+2}}\,.
\end{equation}
It follows that eq.~(\ref{eq:tspectrum}) can be numerically inverted by any root-finding algorithm, provided $0<\bmu<\rho/(v+2)$. We can regard eq.~(\ref{eq:tspectrum}) as a rough approximation to the original problem. In order to use it when $\mu_\text{min} \equiv \min_k\{\mu_k\} < \max_k\{\mu_k\} \equiv \mu_\text{max}$, we must preliminarily define $\tilde\mu$ in terms of the components of~$\mu$. One possibility is to average these, {\it i.e.} to choose
\begin{equation}
\tilde \mu \equiv \frac{1}{v}\sum_{k=1}^v\mu_k\,.
\label{eq:averchoice}
\end{equation}
Subject to this choice, we expect $\tilde\lambda$ to lie somewhere between $\lambda_\text{min} \equiv \min_k\{\lambda_k\}$ and $\lambda_\text{max}\equiv\max_k\{\lambda_k\}$. When $\mu_\text{min}<\mu_\text{max}$, eq.~(\ref{eq:tspectrum}) can be thought of as the lowest order approximation of a perturbative expansion of eqs.~(\ref{eq:truncspectrum}) around the point $\lambdaT = \{\blambda,\ldots,\blambda\}$. If the ratio $\mu_\text{max}/\mu_\text{min}$ is not exceedingly large, such expansion is expected to converge quickly. Few perturbative corrections to $\lambdaT$ should be sufficient to guarantee a good level of approximation. In view of eq.~(\ref{eq:averchoice}), the inverse operator $\cT_\rho^{-1}$ is defined within the convex domain
\begin{equation}
  \cD(\cT_\rho^{-1}) = \left\{\mu:\ \sum_{k=1}^v \mu_k < \frac{\rho v}{v+2}\right\}\,.
\end{equation}
To allow for a comparison, we also plot $\cD(\cT_\rho^{-1})$ for $v=3$ in Fig.~\ref{fig:defdoms} (right).  Aim of the present paper is to carry out a theoretical study of the perturbative expansion of eqs.~(\ref{eq:truncspectrum}) up to the fourth order and to discuss some analytic aspects of it. 

There are several reasons why the suggested perturbative approach can outperform the iterative procedure introduced in ref.~\cite{palombi4}. In first place, it can be easily checked that $\cD(\tau_\rho^{-1})\subset \cD({\cal T}_\rho^{-1})$. Accordingly, perturbative approximations to $\tau_\rho^{-1}$ are well defined even when statistical fluctuations in sample space make $\mu\notin\cD(\tau_\rho^{-1})$, provided at least $\mu\in\cD(\cT_\rho^{-1})$. In this sense, the perturbative approach represents a regularization of the covariance reconstruction problem with respect to the existence of a solution  (see sect.~7 of ref.~\cite{palombi4} for a detailed discussion of the ill-posedness of the problem).

Secondly, statistical fluctuations of the higher components of $\mu$ are always amplified by eqs.~(\ref{eq:truncspectrum}) as a
consequence of the non-linearity and the unboundedness of $\tau_\rho^{-1}$, sometimes resulting in unacceptably large variances for the higher components of $\lambda$. In the literature of inverse problems, this property is known as the {\it instability} of the inverse operator, see for instance
ref.~\cite{cavalier}. In the
framework of perturbation theory, non-linearity arises systematically while increasing the order of the approximation, since each perturbative correction depends non-linearly upon the previous ones. Therefore, by stopping the expansion at different orders we have the possibility to define a class of statistical estimators of $\lambda$, each characterized by its own bias and variance. We find that the variance of the upper (lower) half of the reconstructed eigenvalues increases (decreases) together with the perturbative order, whereas their bias always decreases. In all applications where the upper part of $\lambda$ matters, it is therefore possible ---at least in principle--- to optimize the choice of the perturbative estimator according to one's specific needs. This provides a regularization of the reconstruction problem with respect to the stability properties of the solution. 

Finally, in critical applications requiring reconstructions of the covariance matrix for several values of $\rho$ and/or $\mu$, it could be important to get fast yet approximate solutions, such as perturbation theory provides, rather than slow yet exact ones. Indeed,  the convergence of the fixed point iteration was shown to slow down as  $\rho\to 0$ or $v\to\infty$, the rate of the slowing down being polynomial in  the former limit and exponential in the latter. Thus, in all situations where $\rho\lesssim\lambda_\text{min}$ or $v\gg 1$, the use of the fixed point algorithm could be unfavorable. 

The plan of the paper is as follows. In sect.~2 we derive some preliminary results concerning the integrals $\alpha$ and $\alpha_k$, which are necessary for a systematic implementation of the perturbative strategy to all orders. In sect.~3 we work out the expansion by means of paper-and-pencil calculations and \Maple programs (reported in App.~C). Sect.~4 is devoted to discussing some analytic properties of the first few perturbative terms, when $\bmu$ is chosen as in eq.~(\ref{eq:averchoice}). In sect.~5 we simulate the statistical properties of the perturbative estimators of $\lambda$ and the iterative one in sample space for a specific choice of the covariance matrix. We finally draw our conclusions in sect.~6.    

\section{Building blocks in Tallis' limit}

As well known, perturbation theory is an expansion technique around a reference solution, which is assumed to be either calculable or easily computable. Its mathematical structure develops from building blocks which are themselves entirely defined in terms of the reference solution. When it comes to perturbing eqs.~(\ref{eq:truncspectrum}) around $\lambdaT$, the elementary objects we need to focus on are the Gaussian integrals $\alpha_{k\ell m\ldots}(\rho;\lambda)$ and their partial derivatives for $\lambda\to\lambdaT$. In the sequel we refer to this limit as \emph{Tallis' limit}. 

To begin with, we set-up the notation. Throughout the paper we let $\lambdaT = \{\blambda,\ldots,\blambda\}$ denote a set of fully degenerate eigenvalues. Outside Tallis' limit a convenient axes reshuffling allows us to assume the ordering $\mu_1\le\mu_2\le\ldots\le\mu_v$ with no loss of generality (this also follows naturally from assuming $\lambda_1\le\lambda_2\le\ldots\le\lambda_v$, as asserted by Proposition 5 of ref.~\cite{palombi4}). We denote a generic Gaussian integral with $n$ (not necessarily distinct) indices $\{i_1,\ldots,i_n\}\equiv {\cal I}$ either by the notation $\alpha_{i_1\ldots i_n}(\rho;\lambda)$ introduced in sect.~1 or by the {\it compact} notation $\alpha_{1:m_1\ldots v:m_v}(\rho;\lambda)$, where the generic subscript $k\!:\!m_k$ indicates that the directional index $k$ occurs with multiplicity $m_k$ in $\cI$. Similarly, we denote the $n^\text{th}$~order derivative operator with respect to the variances $\lambda_{i_1}\ldots,\lambda_{i_n}$ either by the {\it standard} symbol $\partial_{i_1\ldots i_n} = \partial^n/(\partial  \lambda_{i_1}\ldots \partial\lambda_{i_n})$ or by its {\it compact} version $\partial_{1:m_1\ldots v:m_v} = {\partial^{\,m_1+\ldots+m_v}}/{(\partial\lambda_1)^{m_1}\ldots(\partial\lambda_v)^{m_v}}$. Given $\cI$, the multiplicity set $\cM_\cI = \{m_1,\ldots,m_v\}$ in unambiguously defined and the other way round. For consistency $\cM_\cI$ fulfills the constraint $\sum_{k=1}^v m_k = n$. In case of vanishing multiplicities we omit to write all the corresponding indices. For instance, we write $\alpha_{k:m_k}(\rho;\lambda)$ in place of the rather pedantic $\alpha_{1:0\ldots k:m_k\ldots v:0}(\rho;\lambda)$ as well as $\partial_{k:m_k}$ in place of $\partial_{1:0\ldots k:m_k\ldots v:0}$. Last but not least, we refer the reader to ref.~\cite{palombi4} for definitions and properties concerning the truncation operator $\tau_\rho$ and its inverse~$\tau_\rho^{-1}$.

Having said that, we start our investigation of Tallis' limit with a simple proposition, meant to sum up the results of ref.~\cite{tallis}:

\begin{prop}
For $m_1\ge 0$, \ldots, $m_v\ge 0$, we have
\begin{equation}
\label{eq:alphalimone}
\alpha_{1:m_1\ldots v:m_v}(\rho;\lambdaT) = \Delta_{1:m_1\ldots v:m_v}\,
F_{v+2(m_1+\ldots+m_v)}\left(\frac{\rho}{\blambda}\right)\,,
\end{equation}
with
\begin{equation}
\Delta_{1:m_1\ldots v:m_v} = \prod_{j=1}^v(2m_j-1)!!\,.
\label{eq:deltasymbol}
\end{equation}
\end{prop}
\begin{proof}
In order to derive eq.~(\ref{eq:alphalimone}), we represent $\alpha_{1:m_1\ldots v:m_v}(\rho;\lambdaT)$ in spherical coordinates,   
\begin{align}
x_1 & = rf_1(\theta_1,\ldots,\theta_{v-1})\,,\nonumber\\[-1.0ex]
& \ \, \vdots \label{eq:sphercoords}\\[-1.0ex]
x_v & = rf_v(\theta_1,\ldots,\theta_{v-1})\,.\nonumber
\end{align}
Recall that $\sum_k f_k^2 = 1$ and $\rd^v x = r^{v-1}\rd r\rd\Omega$, with $\rd\Omega$ embodying the angular part of the Jacobian of eq.~(\ref{eq:sphercoords}) and the differentials of the angles $\theta_1$, \ldots, $\theta_{v-1}$. With a bit of algebra we easily arrive at
\begin{equation}
\alpha_{1:m_1\ldots v:m_v}(\rho;\lambdaT) = \Delta_{1:m_1\ldots v:m_v}F_{v+2(m_1+\ldots+m_v)}\left(\frac{\rho}{\blambda}\right)\,,
\end{equation}
where $\Delta_{1:m_1\ldots v:m_v}$ is just a proportionality factor independent of $\rho$ and $\blambda$. In order to fix $\Delta_{1:m_1\ldots v:m_v}$, we observe that $\alpha_{1:m_1\ldots v:m_v}(\rho;\lambdaT)$ factorizes into $v$ one-dimensional integrals as $\rho\to\infty$, corresponding to unconditioned univariate Gaussian moments of orders $2m_1$,\ldots,$2m_v$, normalized respectively by powers $\blambda^{m_1}$,\ldots,$\blambda^{m_v}$ of
the common variance. Hence, we conclude that $\Delta_{1:m_1\ldots v:m_v} = \prod_{k=1}^v(2m_k-1)!!$.  
\end{proof}

The $\smash{\chi^2}$-\cdf $\smash{F_{v+2(m_1+\ldots+m_v)}(\rho/\blambda)}$ is intuitively interpreted as a correction factor incorporating all the effects of integrating Gaussian densities on a finite volume. Obviously, it lessens the value of the integral except for $\rho\to\infty$, a limit for which it converges to one.

Whenever Gaussian integrals are expressed in {\it standard} notation, eq.~(\ref{eq:alphalimone}) can be applied only provided the index multiplicities are preliminarily counted. Nevertheless, it is reasonable to expect $\smash{\Delta_{1:m_1\ldots v:m_v}}$ to be the {\it compact} representation of some not yet specified coefficient $\smash{\Delta_{i_1\ldots i_n}}$. By the same argument used above, i.e. by letting $\rho\to\infty$, we conclude that the latter is given by 
\begin{equation}
 \Delta_{i_1\ldots i_n} = \E[z_{i_1}^2\ldots z_{i_n}^2]\,,\quad\text{with}\quad
 z_{i_k}\, \stackrel{\footnotesize \rm iid}{\sim}\,\cN(0,1)\,,\qquad k=1,\ldots,n\,.
\label{eq:isserlis}
\end{equation}
Thanks to Isserlis' Theorem\footnote{In mathematical physics the same result is universally known as {\it Wick's Theorem}.}~\cite{isserlis}, $\Delta_{i_1\ldots i_n}$ can be expressed as a sum of products of Kronecker symbols, namely
\begin{equation}
\Delta_{i_1\ldots i_n} = \sum\prod \delta_{ij}\,,
\label{eq:DeltaKron}
\end{equation}
where the sum includes all distinct ways of partitioning the duplicated set $\cI^2\equiv \{i_1,i_1,\ldots,i_n,i_n\}$ into pairs and the product is over all pairs of a single partition. For later convenience it is worthwhile listing the lowest order coefficients:
\begin{align}
 \Delta_{i_1} & = 1\,, \\[1.0ex]
\Delta_{i_1i_2} & = 1 + 2\delta_{i_1i_2}\,, \\[1.0ex]
 \Delta_{i_1i_2i_3} & = 1 + 2\left(\delta_{i_1i_2}\!+\!
  \delta_{i_1i_3}\!+\!\delta_{i_2i_3}\right)\!+\!2\left(\delta_{i_1i_2}\delta_{i_1i_3}\!+\!\delta_{i_1i_2}\delta_{i_2i_3}
\!+ \!\delta_{i_1i_3}\delta_{i_1i_2}\right)\!
 +\! 2\,\delta_{i_1i_2}\delta_{i_1i_3}\delta_{i_2i_3}\,,\\[0.0ex]
& \ \, \vdots\nonumber
\end{align}
Eqs.~(\ref{eq:isserlis})--(\ref{eq:DeltaKron}) allow us to reformulate
Proposition 2.1 as follows:

\vskip\parskip

\begin{corol} For $n\ge 0$ and $\{i_1,\ldots,i_n\}$ a set of $n$ (not necessarily distinct) indices, we have
\begin{equation} 
\alpha_{i_1\ldots i_n}(\rho;\lambdaT) = \Delta_{i_1\ldots i_n}F_{v+2n}\left(\frac{\rho}{\blambda}\right)\,,
\label{eq:alphalimtwo}
\end{equation}
with $\Delta_{i_1\ldots i_n}$ as in eq.~(\ref{eq:DeltaKron}).
\end{corol}

Now, in order to take arbitrarily high order derivatives of $\alpha_{i_1\ldots i_n}(\rho;\lambda)$ with respect to the components of $\lambda$, we iterate the basic rule
\begin{equation}
(2\lambda_k\partial_k)\,\alpha_{i_1\ldots i_n}(\rho;\lambda) = \alpha_{i_1\ldots
  i_n k}(\rho;\lambda) - \biggl(1+2\sum_{j=1}^n\delta_{ki_j}\biggr)\alpha_{i_1\ldots i_n}(\rho;\lambda)\,,
\label{eq:basicder}
\end{equation}
which follows from differentiating $\alpha_{i_1\ldots i_n}$ under the integral sign. Note that eq.~(\ref{eq:basicder}) is not specifically related to spherical truncations, {\it i.e.} it is formally invariant under a reshaping of the truncation surface. It also shows that the differential operator $2\lambda_k\partial_k$ behaves in a simpler manner than $\partial_k$ when acting on $\alpha_{i_1\ldots i_n}$, in that it produces an integer linear combination  of Gaussian integrals. For this reason, the recurrence generated by $\partial_{k:n}$ acting on $\alpha_{i_1\ldots i_n}$ can be derived by first working out the action of the operator $(2\lambda_k\partial_k)^n$ on $\alpha_{i_1\ldots i_n}$ and by then using
\begin{equation}
  \partial_{k:n} = \frac{1}{\lambda_k^n}\sum_{j=1}^n \frac{(-1)^{n-j}}{2^j} {n\brack j} (2\lambda_k\partial_k)^j\,,
  \label{eq:oprelation}
\end{equation}
where the symbols ${n\brack j}$ represent unsigned Stirling numbers of the first kind. Eq.~(\ref{eq:oprelation}) is a standard result of combinatorial analysis. We refer the reader to exercise 13, chap.~6 of ref.~\cite{knuth} for a proof of it. Iterated applications of the operator $2\lambda_k\partial_k$ generate increasingly involved sums of Gaussian integrals, as asserted by 

\begin{prop} For all $j\ge 1$ and $n\ge 0$, we have
\begin{equation}
(2\lambda_k\partial_k)^j\alpha_{k:n} = \sum_{r=0}^j(-1)^{j-r}c_{jr}(n)\alpha_{k:(n+r)}\,,
\label{eq:partderone}
\end{equation}
with
\begin{equation}
c_{jr}(n) = \sum_{\ell_1=0}^r\sum_{\ell_2=0}^{\ell_1}\,\ldots
\sum_{\ell_{j-r}=0}^{\ell_{j-r-1}}\ \prod_{s=1}^{j-r}[2(n+\ell_s)+1]\,.
\label{eq:expcoef}
\end{equation}
\end{prop}
\begin{proof}
The proof is by induction. We first note that $c_{10}(n) = (2n+1)$ and  $c_{11}(n)=1$. Hence, for $j=1$ eq.~(\ref{eq:partderone}) agrees with
eq.~(\ref{eq:basicder}). Now, suppose that $(2\lambda_k\partial_k)^j\alpha_{k:n}$ is well represented by eq.~(\ref{eq:partderone}) with $c_{jr}(n)$ as in eq.~(\ref{eq:expcoef}). Then, 
\begin{align}
(2\lambda_k\partial_k)^{j+1}\alpha_{k:n} & =
\sum_{r=0}^j(-1)^{j-r}c_{jr}(n)(2\lambda_k\partial_k)\alpha_{k:(n+r)} \nonumber\\[0.0ex]
& = \sum_{r=0}^j(-1)^{j-r}c_{jr}(n)\{\alpha_{k:(n+r+1)} -
[2(n+r)+1]\alpha_{k:(n+r)}\} \nonumber\\[0.0ex]
& = \sum_{r=0}^{j+1}(-1)^{j+1-r}\{c_{j(r-1)}(n) + [2(n+r)+1]c_{jr}(n)\}\alpha_{k:(n+r)}\,.
\end{align}
Hence, the proof is complete if we are able to show that $c_{jr}(n)$ fulfills the recurrence
\begin{equation}
c_{(j+1)r}(n) = c_{j(r-1)}(n) + [2(n+r)+1]c_{jr}(n)\,,
\label{eq:crecur}
\end{equation}
To this aim, we first calculate the second term on the \rhs as
\begin{align} 
[2(n+r)+1]c_{jr}(n)
& =
\sum_{\ell_1=0}^{r}\sum_{\ell_2=0}^{\ell_1}\ldots\!\!\sum_{\ell_{j-r}=0}^{\ell_{j-r-1}}[2(n+r)+1)]\prod_{s=1}^{j-r}[2(n+\ell_s)+1]\nonumber\\[2.0ex]
& = \sum_{\ell_2=0}^{r}\sum_{\ell_3=0}^{\ell_2}\ldots\!\!\!\!\sum_{\ell_{j+1-r}=0}^{\ell_{j-r}}[2(n+r)+1)]\prod_{s=2}^{j+1-r}[2(n+\ell_s)+1]\nonumber\\[2.0ex]
& = \sum_{\ell_1=r}^r\sum_{\ell_2=0}^{\ell_1}\ldots\!\!\!\!\sum_{\ell_{j+1-r}=0}^{\ell_{j-r}}\prod_{s=1}^{j+1-r}[2(n+\ell_s)+1]\,,
\end{align}
and then we add it to the first one, thus obtaining 
\begin{align}
c_{j(r-1)}(n) & + [2(n+r)+1]c_{jr}(n) =
\sum_{\ell_1=0}^{r-1}\sum_{\ell_2=0}^{\ell_1}\ldots\!\!\!\!\sum_{\ell_{j+1-r}=0}^{\ell_{j-r}}\prod_{s=1}^{j+1-r}[2(n+\ell_s)+1]\nonumber\\[1.0ex]
&  +
\sum_{\ell_1=r}^r\sum_{\ell_2=0}^{\ell_1}\ldots\!\!\!\!\sum_{\ell_{j+1-r}=0}^{\ell_{j-r}}\prod_{s=1}^{j+1-r}[2(n+\ell_s)+1]
= c_{(j+1)r}(n)\,.
\end{align}
\end{proof}
Nested sums similar to eq.~(\ref{eq:expcoef}) are considered for instance in ref.~\cite{butler}, where all cases in study are reduced to closed-form expressions with the help of special numbers, such as binomial coefficients, Stirling numbers, center factorial numbers, {\it etc}. Eq.~(\ref{eq:expcoef}) looks a bit harder to manage, since the summand is a product of non-homogeneous functions of the sum variables, hence it is not clear whether the nested sums can be ultimately evaluated in closed-form. However, as far as we are concerned, a convenient representation of $c_{jr}(n)$ is provided by  

\begin{prop} For all $j\ge 1$, $0\le r\le j$ and $n\ge 0$, we have
\begin{equation}
c_{jr}(n) = \sum_{t=r}^j(-2)^{j-t}\frac{[2(n+t)-1]!!}{[2(n+r)-1]!!}{j\brace t}{t
  \choose r}\,,
\label{eq:resum}
\end{equation}
with the symbols ${j\brace t}$ denoting Stirling numbers of the second kind.
\end{prop}
\begin{proof}
We let $d_{jr}(n)$ denote the \rhs of eq.~(\ref{eq:resum}). For $j=1$, we have $d_{10}(n) = 2n+1$ and $d_{11}(n) = 1$, which equal respectively $c_{10}(n)$ and $c_{11}(n)$. Thus, we just need to prove that $d_{jr}(n)$ fulfills eq.~(\ref{eq:crecur}). To this end, it is sufficient to make use of the basic recursive formulae ${n+1 \brace m} = m{n\brace m} + {n\brace m-1}$ and ${n+1 \choose m} = {n\choose m} + {n\choose m-1}$. We detail the algebra for the sake of completeness. We start from
\begin{align}
d_{(j+1)r}(n) =
\sum_{t=r}^{j+1}(-2)^{j+1-t}\frac{[2(n+t)-1]!!}{[2(n+r)-1]!!}{j+1 \brace
  t}{t\choose r} \nonumber
\end{align}
\begin{align}
\hskip 1.5cm & = \sum_{t=r}^{j}(-2)^{j-t}\frac{[2(n+t)-1]!!}{[2(n+r)-1]!!}(-2t){j \brace
  t}{t\choose r} \nonumber\\[2.0ex]
& + \sum_{t=r-1}^{j}(-2)^{j-t}\frac{[2(n+t)+1]!!}{[2(n+r)-1]!!}{j \brace
  t}{t+1\choose r}\,\nonumber\\[2.0ex]
& = (2n+1)d_{jr}(n) + \sum_{t=r-1}^{j}(-2)^{j-t}\frac{[2(n+t)+1]!!}{[2(n+r)-1]!!}{j \brace
  t}{t\choose r-1}\,.
\label{eq:drecrel}
\end{align}
Then, we add and subtract $(2r)d_{jr}(n)$ to the \rhs of
eq.~(\ref{eq:drecrel}). Hence, 
\begin{align}
d_{(j+1)r}(n) & =[2(n+r)+1]d_{jr}(n)\nonumber\\[2.0ex]
& + \sum_{t=r-1}^{j}(-2)^{j-t}\frac{[2(n+t)-1]!!}{[2(n+r)-1]!!}{j \brace
  t}\left[[2(n+t)+1]{t\choose r-1}-2r{t\choose r}\right]\,.
\end{align}
Since $r{t\choose r} = (t-r+1){t\choose r-1}$, the second term on the \rhs is recognized to be $d_{j(r-1)}(n)$. 
\end{proof}

In view of ref.~\cite{butler} it is no surprise that $c_{jr}(n)$ can be represented in terms of Stirling numbers. In addition, the presence of terms such as $j\brace t$ fits perfectly when combining eq.~(\ref{eq:oprelation}) with eq.~(\ref{eq:partderone}). We recall indeed that Stirling numbers of the first and second kind fulfill the identity
\begin{equation}
\sum_{t=0}^{\max\{j,k\}}(-1)^{t-k}{t\brace j}{k \brack t} = \delta_{jk}\,.
\label{eq:stirlinv}
\end{equation}

We have now collected all the ingredients needed to prove the main result of this section, namely

\noindent\colorbox{shadecolor}
{
  \parbox{1.0\textwidth}
  {
    \begin{prop} For $m,n\ge 0$, let $\cK = \{k_1,\ldots,k_m\}$ and $\cI = \{i_1,\ldots,i_n\}$ denote sets of (not
      necessarily distinct) directional indices, {\it i.e.} $1\le k_j\le v$ and $1\le i_j\le v$. We have
      \begin{equation}
        \label{eq:dergeneral}
        \partial_{k_1...k_m}\alpha_{i_1\ldots i_n}(\rho;\lambdaT) = \frac{1}{2^m\blambda^m}\Delta_{k_1\ldots
          k_mi_1\ldots i_n}\sum_{j=0}^{m}(-1)^{m-j}\binom{m}{j}F_{v+2(j+n)}\left(\frac{\rho}{\blambda}\right)\,.
      \end{equation}
      In particular, for $n=0,1$ we have
      \begin{align}
        \label{eq:deralpha}
        \partial_{k_1...k_m}\alpha(\rho;\lambdaT) & = \frac{1}{2^m\blambda^m}\Delta_{k_1\ldots
          k_m}\sum_{j=0}^{m}(-1)^{m-j}\binom{m}{j}F_{v+2j}\left(\frac{\rho}{\blambda}\right)\,,\\[1.0ex]
        \label{eq:deralphak}
        \partial_{k_1...k_m}\alpha_i(\rho;\lambdaT) & = \frac{1}{2^m\blambda^m}\Delta_{k_1\ldots k_m
          i}\sum_{j=0}^{m}(-1)^{m-j}\binom{m}{j}F_{v+2(j+1)}\left(\frac{\rho}{\blambda}\right)\,.
      \end{align}
    \end{prop}
  }
}

\begin{proof} We let $\cM_\cK = \{m_1,\ldots,m_v\}$ and $\cM_\cI = \{n_1,\ldots,n_v\}$ denote the multiplicity sets associated respectively to $\cK$ and $\cI$, such that $\partial_{k_1\ldots k_m}=\partial_{1:m_1\ldots v:m_v}$,  $\alpha_{i_1\ldots i_n} = \alpha_{1:n_1\ldots v:n_v}$ and
\begin{equation}
\partial_{k_1...k_m}\alpha_{i_1\ldots i_n}(\rho;\lambdaT) = \partial_{1:m_1...v:m_v}\alpha_{1:n_1\ldots v:n_v}(\rho;\lambdaT)\,.
\end{equation}
We recall that $\cM_\cK$ and $\cM_\cI$ fulfill $\sum_{\ell=1}^v m_\ell = m$ and $\sum_{\ell=1}^v n_\ell = n$ for consistency. Using eq.~(\ref{eq:oprelation}) $m$ times yields 
\begin{align}
\partial_{k_1...k_m}\alpha_{i_1\ldots i_n}(\rho;\lambdaT) & = 
\frac{1}{\blambda^m}\sum_{j_1=1}^{m_1}\frac{(-1)^{m_1-j_1}}{2^{j_1}}{m_1\brack
  j_1}\,\ldots\,\sum_{j_v=1}^{m_v}\frac{(-1)^{m_v-j_v}}{2^{j_v}}{m_v\brack
  j_v}\nonumber\\[1.0ex]
& \ \cdot\ \left[(2\lambda_1\partial_1)^{j_1}\,\ldots\,
  (2\lambda_v\partial_v)^{j_v}\alpha_{1:n_1\ldots v:n_v}\right](\rho;\lambdaT)\,.
\label{eq:dummopsum}
\end{align}
Moreover, with the help of eqs.~(\ref{eq:partderone}) and (\ref{eq:alphalimone}), eq.~(\ref{eq:dummopsum})
reduces to
\begin{align}
& \partial_{k_1...k_m}\alpha_{i_1\ldots i_n}(\rho;\lambdaT) \nonumber\\[1.5ex]
& \hskip 0.7cm =  \frac{1}{\blambda^m}\sum_{j_1=1}^{m_1}\frac{1}{2^{j_1}}{m_1\brack
  j_1}\sum_{q_1=0}^{j_1}{(-1)^{m_1-q_1}}\ldots \sum_{j_v=1}^{m_v}\frac{1}{2^{j_v}}{m_v\brack
  j_v}\sum_{q_v=0}^{j_v}{(-1)^{m_v-q_v}}\nonumber\\[1.0ex]
& \hskip 0.7cm \ \cdot \ [2(n_1+q_1)-1]!!\,\ldots\,[2(n_v+q_v)-1]!! \nonumber\\[1.0ex]
& \hskip 0.7cm \ \cdot\, c_{j_1q_1}(n_1)\, \ldots\, c_{j_vq_v}(n_v) \,\cdot\, F_{v+2(s+q_1+\ldots+q_v)}\left(\frac{\rho}{\blambda}\right)\,. 
\end{align}
As a next step, we evaluate the coefficients $\smash{c_{j_1q_1}(n_1)}$, \ldots, $\smash{c_{j_vq_v}(n_v)}$ in terms of the expressions obtained in Proposition 2.3. Then we make use of eq.~(\ref{eq:stirlinv}) to make all Stirling numbers disappear. We finally identify $\smash{\prod_{j=1}^v[2(n_j+m_j)-1)]!! = \Delta_{1:(m_1+n_1)\ldots v:(m_v+n_v)}}$ and so we arrive at  
\begin{align}
& \partial_{k_1...k_m}\alpha_{i_1\ldots i_n}(\rho;\lambdaT) = \frac{1}{2^m\blambda^m}\Delta_{1:(m_1+n_1)\ldots
    v:(m_v+n_v)}\nonumber\\[2.0ex]
& \hskip 0.7cm \cdot \, \sum_{q_1=0}^{m_1}\ldots\sum_{q_v=0}^{m_v}
(-1)^{m-(q_1+\ldots + q_v)}{m_1 \choose q_1}\ldots {m_v \choose q_v}F_{v+2(s+q_1+\ldots+q_v)}\left(\frac{\rho}{\blambda}\right)\,.
\label{eq:almfin}
\end{align}
Two additional steps are needed to complete the proof. In first place, we notice that the multiplicity set associated to $\cK\cup\cI = \{k_1,\ldots,k_m,i_1,\ldots,i_n\}$ is $\cM_{\cK\cup\cI} = \{m_1+n_1,\ldots,m_v+n_v\}$. Therefore, by Isserlis' Theorem we identify $\Delta_{1:(n_1+m_1)\ldots  v:(n_v+m_v)} = \Delta_{k_1\ldots k_mi_1\ldots i_n}$. In second place, we observe that $(q_1+\ldots+q_v)$ ranges from 0 to $m$ for $0\le q_1\le m_1$, \ldots, $0\le q_v\le m_v$. Hence, we recast the \rhs of eq.~(\ref{eq:almfin}) to 
\begin{equation}
\partial_{k_1...k_m}\alpha_{i_1\ldots i_n}(\rho;\lambdaT) =
\frac{1}{2^m\blambda^m} \Delta_{k_1\ldots k_mi_1\ldots i_n} \,
\sum_{j=0}^{m}(-1)^{m-j} e_j
F_{v+2(s+j)}\left(\frac{\rho}{\blambda}\right)\,, 
\end{equation}
with
\begin{equation}
e_j = \sum_{q_1=0}^{m_1}\ldots\sum_{q_v=0}^{m_v} {m_1 \choose q_1}\ldots {m_v
  \choose q_v} \delta_{q_1+\ldots+q_v,j}\,.
\end{equation}
This multiple sum is easily calculated by iterating the Vandermonde's convolution $\sum_k {r \choose k}{s\choose p-k} = {r+s \choose p}$. This finally yields $e_j = {m\choose j}$.  
\end{proof}

\section{Perturbative expansion}

In order to solve eqs.~(\ref{eq:truncspectrum}), perturbation theory prescribes that we interpret $\mu$ and $\lambda$ as smooth functions of a
parameter~$\epsilon\in [0,1]$, such that $\lambda(\epsilon=0) = \lambdaT$ and $\lambda(\epsilon=1) = \tau_{\rho}^{-1}\cdot \mu$. We must consider
$\epsilon$ as an auxiliary variable, allowing us to pass continuously from Tallis' limit to the ultimate solution we are seeking. Then, we are requested to expand $\mu(\epsilon)$ and $\lambda(\epsilon)$ in power series of $\epsilon$ around the point $\epsilon=0$, namely
\begin{alignat}{8}
\label{eq:lambdaexp}
\lambda_k(\epsilon)\, & =\, \blambda\, & + &\, \epsilon \lambda_k^{(1)} & + &\, \epsilon^2
\lambda_k^{(2)} & + & \ldots\,,\\[2.0ex]
\label{eq:muexp}
\mu_k(\epsilon)\, & =\, \bmu\, & + &\, \epsilon \mu_k^{(1)} & + &\, \epsilon^2 \mu_k^{(2)} &
+ & \ldots\,.
\end{alignat}
For later convenience we let $R_k$ denote the integral ratio $\alpha_k/\alpha$. Since $R_k$ depends smoothly upon $\lambda$, it can be analogously expanded in power series of $\epsilon$. Thus, eqs.~(\ref{eq:truncspectrum}) read
\begin{equation}
\bmu+ \epsilon \mu_k^{(1)} + \epsilon^2 \mu_k^{(2)} + \ldots\ = \left(
  \blambda + \epsilon \lambda_k^{(1)} + \epsilon^2\lambda_k^{(2)} +
  \ldots\right)\cdot\left(\bR+ \epsilon R_k^{(1)} + \epsilon^2 R_k^{(2)} +
  \ldots\right)\,,
\label{eq:spectrumexp}
\end{equation}
with $\bR = R_k(\rho;\lambdaT)$ and $R_k^{(n)} =
(n!)^{-1}{\rd^nR_k}/{\rd \epsilon^n}|_{\epsilon=0}$ for $n=1,2,\ldots$\  The idea
underlying perturbation theory is that we treat separately terms in eq.~(\ref{eq:spectrumexp}) belonging to different perturbative orders, that is to say we equal terms of the same order in $\epsilon$ on both sides of eq.~(\ref{eq:spectrumexp}) and then we solve one by one the algebraic equations thus obtained. A few caveats must be noticed.

\noindent {\textbf{\emph{i})} Since $\mu_k$ is an input parameter, we must specify how it enters the Taylor coefficients of the perturbative function $\smash{\mu_k(\epsilon)}$. In principle, the assignment can be made in complete freedom. For instance, the choice we adopt in the sequel
is to confine $\mu_k$ to the first-order Taylor coefficient, namely
\begin{equation}
\left\{\begin{array}{ll}
\mu_k^{(0)} & \!\!\! = \ \bmu\,,\\[1.0ex]
\mu_k^{(1)} & \!\!\! = \ \mu_k - \bmu \, \equiv\, \delta\mu_k\,,\\[1.2ex]
\mu_k^{(n)} & \!\!\! = \ 0\,,\qquad n= 2,3,\ldots
\label{eq:splittingone}
\end{array}\right.
\end{equation}
Alternatively, we might spread $\mu_k$ over all Taylor coefficients of $\mu_k(\epsilon)$, by letting for instance
\begin{equation}
\left\{\begin{array}{ll}
\mu_k^{(0)} & \!\!\! = \ \bmu\,,\\[1.0ex]
\mu_k^{(n)} & \!\!\! = \dfrac{1}{n!}\{\log[1+(\mu_k-\bmu)]\}^n\,,\qquad n = 1,2,\ldots
\end{array}\right.
\label{eq:splittingtwo}
\end{equation}
Both eqs.~(\ref{eq:splittingone}) and (\ref{eq:splittingtwo}) comply with the requirement $\mu_k(\epsilon=1) = \mu_k$. It must be noticed, however, that each legitimate splitting of $\mu$ affects differently the convergence properties of the perturbative series of $\lambda$ as well as the statistical properties of the truncated series, when $\mu$ is turned into a random variable in sample space. This will be further investigated in sect.~5. 

\noindent {\textbf{\emph{ii})} The specific choice of $\bmu$ is rather arbitrary: as far as we are concerned with the feasibility of the perturbative expansion, the only requirement to fulfill is that eq.~(\ref{eq:tspectrum}) be invertible, which is guaranteed provided $0\le\bmu<\rho/(v+2)$. We mentioned in sect.~1 that a convenient choice is represented by     
\begin{equation}
  \bmu = \frac{1}{v}\sum_{k=1}^v \mu_k \equiv \bar\mu\,.
\label{eq:muaver}
\end{equation}
We shall see shortly what benefits derive from eq.~(\ref{eq:muaver}). Certainly, we have $\bar\mu<\rho/(v+2)$ whenever $\mu\in\cD(\tau_\rho^{-1})$, hence eq.~(\ref{eq:muaver}) is at least a legitimate choice. Indeed, if $\mu\in\cD(\tau_\rho^{-1})$ then from eq.~(\ref{eq:domtauromo}) it follows that
\begin{equation}
  \sum_{k=1}^{v}\left(\sum_{j\ne k}^{1\ldots v}\mu_j + 3\mu_k\right) = v(v+2)\bar\mu < \sum_{k=1}^v\rho = v\rho \qquad \Rightarrow \qquad \bar\mu< \rho/(v+2)\,.
\end{equation}

\noindent {\textbf{\emph{iii})} Perturbation theory works only provided the $\Oen$-equations
\begin{equation}
\cE_k^{(n)}\equiv\sum_{j=0}^n \lambda_k^{(n-j)}R_k^{(j)} - \mu_k^{(n)} = 0\,,\qquad n=0,1,\ldots
\label{eq:Oneq}
\end{equation}
obtained by collecting all the $\Oen$-terms from  eq.~(\ref{eq:spectrumexp}), yield an algebraic relation among the Taylor coefficients of
$\lambda(\epsilon)$ which can be solved with respect to $\lambda^{(n)}$. This allows us to represent the latter as a function $\lambda^{(n)}(\blambda,\lambda^{(1)},\ldots,\lambda^{(n-1)},\mu^{(n)})$ of the lower order coefficients of $\lambda(\epsilon)$ together with~$\mu^{(n)}$. If such property is fulfilled, as we argue in a while, then solving eqs.~(\ref{eq:truncspectrum}) perturbatively means solving one after another the systems of equations
\begin{equation}
  \{{\cal E}_k^{(1)}=0\}_{k=1}^v\,, \quad \{{\cal E}_k^{(2)}=0\}_{k=1}^v\,, \quad \ldots \quad \,, \quad \{{\cal E}_k^{(n)}=0\}_{k=1}^v\,,
\end{equation}
up to a predefined order $n$. Establishing the level of precision thus achieved is a complicated task, as typical of asymptotic expansions. From a qualitative point of view, the approximation is certainly correct up to $\Oenpo$-terms. However, that is not a quantitative estimate of the truncation error.   

\subsection{General structure of the perturbative expansion}

Regarding point {\bf \emph{iii})}, we notice that the only contributions to $\smash{\cE_k^{(n)}}$ depending explicitly upon $\smash{\lambda^{(n)}}$ are $\smash{\lambda^{(n)}_k\bR}$ and $\smash{\blambda R_k^{(n)}}$. All other terms are of the form $\smash{\lambda_k^{(n-j)}R_k^{(j)}}$ for $j=1,\ldots,n-1$. These terms depend upon $\blambda$, $\lambda^{(1)}$, \ldots, $\lambda^{(n-1)}$, but not upon $\lambda^{(n)}$. Indeed, $R_k$ depends upon $\epsilon$  implicitly via $\lambda(\epsilon)$, thus its $j^{\rm th}$ order derivative with respect to $\epsilon$ distributes progressively according to the chain rule of differentiation. When evaluating the derivative at $\epsilon=0$, all terms proportional to strictly positive powers of $\epsilon$ vanish. As a consequence, each surviving term is proportional to a product of Taylor coefficients of $\lambda(\epsilon)$, each belonging to $\{\blambda,\lambda^{(1)},\ldots,\lambda^{(j)}\}$. In particular, when $j=n$ an explicit calculation yields   
\begin{equation}
\text{terms proportional to $\lambda^{(n)}$ in ${\cal E}_k^{(n)}$}\ :\qquad \lambda^{(n)}_k
\bR + \blambda\sum_{j=1}^v\lambda_j^{(n)}\partial_jR_k(\rho;\lambdaT)\,.
\label{eq:termsln}
\end{equation}
To evaluate the first order partial derivatives of $R_k$, we make use of 
Propositions 2.1 and~2.4. We have
\begin{align}
\partial_jR_k(\rho;\lambdaT) & = \left[\frac{\partial_j\alpha_k}{\alpha} -
  \frac{\alpha_k\partial_j\alpha}{\alpha^2}\right](\rho;\lambdaT)
\nonumber\\[1.0ex]
& = \frac{1}{2\blambda}\left[(1 +
2\delta_{jk})\frac{F_{v+4}}{F_v} - \frac{F_{v+2}^2}{F_v^2} -
2\delta_{jk}\frac{F_{v+2}}{F_v}\right]\biggl(\frac{\rho}{\blambda}\biggr)\,,
\label{eq:firstRder}
\end{align}
hence eqs.~(\ref{eq:Oneq}) reduce to\footnote{From now on we omit to write the argument of the $\chi^2$-\cdf\!'s, which is always understood to be $\rho/\blambda$.}

\noindent\colorbox{shadecolor}
{
  \parbox{0.975\textwidth}
  {
    \vskip-0.3cm
    \begin{equation}
      \sum_{j=1}^v \cJ_{kj}\lambda_j^{(n)} = {\cal G}_k^{(n)}\left(\blambda,\lambda^{(1)},\ldots,\lambda^{(n-1)},\mu^{(n)}\right)\,,
      \label{eq:finalpt}
    \end{equation}
    \\[-6.0ex]
  }
}
\\[1.0ex]
\noindent with
\begin{equation}
\cJ_{kj} = \frac{1}{2}\left[(1 +
2\delta_{kj})\frac{F_{v+4}}{F_v} - \frac{F_{v+2}^2}{F_v^2}\right]\,.
\end{equation}
\begin{center}
\begin{figure}[!t]
\begin{center}
\includegraphics[width=0.7\textwidth]{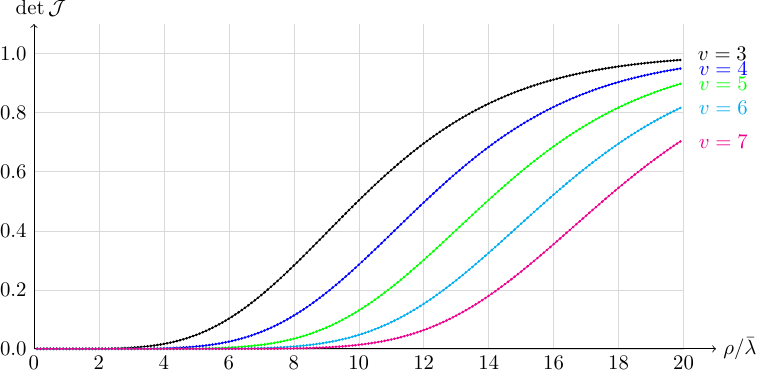}
\vskip 0.2cm
\caption{\small $\det\cJ$ vs. $\rho/\blambda$ at $v=3,\ldots,7$. \label{fig:detJ}}    
\end{center}
\end{figure}
\end{center}
\vskip -0.4cm
We have obtained a system of linear equations with $\lambda^{(n)}$ and
$\cG^{(n)}$ representing respectively the unknown vector and the (known) constant. Moreover, the coefficient matrix $\cJ$ is the Jacobian of the truncation
operator $\tau_\rho$ in Tallis'~limit. Its determinant is given by
\begin{equation}
\det \cJ = \left(\frac{F_{v+4}}{F_v}\right)^{v-1}\left\{
  \left(\frac{v}{2}+1\right)\frac{F_{v+4}}{F_v}-\frac{v}{2}\frac{F_{v+2}^2}{F_v^2}\right\}
\ > \ \frac{2}{v+4}\left(\frac{F_{v+4}}{F_v}\right)^{v-1}\frac{F_{v+2}^2}{F_v^2} > 0 \,,
\label{eq:Jlowerbound}
\end{equation}
the lower bound in eq.~(\ref{eq:Jlowerbound}) following from the inequality $\smash{F_{v+4}F_v/F_{v+2}^2 > (v+2)/(v+4)}$, first proved in ref.~\cite{merkle}. We conclude that $\cJ$ is non-singular for any finite value of $\smash{\rho/\blambda}$, and therefore eq.~(\ref{eq:finalpt}) is unambiguously solved by $\smash{\lambda^{(n)} = \cJ^{-1}\cG^{(n)}}$. Note as well that $\smash{\lim_{\rho/\blambda\to 0}\det\cJ =~0}$, hence the invertibility of $\cJ$ becomes critical at small values of $\smash{\rho/\blambda}$. By way of example, we show in Fig.~\ref{fig:detJ} a plot of $\det \cJ$ vs. $\rho/\blambda$ for $v=3,\ldots,7$. Finally, the inverse of $\cJ$ can be easily worked out. We have indeed
\begin{equation}
(\cJ^{-1})_{jk} = \frac{F_v}{F_{v+4}}\dfrac{
  [(v+2)\delta_{jk}-1]\dfrac{F_{v+4}}{F_v} -
  [v\delta_{jk}-1]\dfrac{F_{v+2}^2}{F_v^2}}{(v+2)\dfrac{F_{v+4}}{F_v} - v\dfrac{F_{v+2}^2}{F_v^2}}\,.
\end{equation}

We now discuss the analytic structure of the known term $\cG_k^{(n)}$. We have just explained that owing to the chain rule of differentiation, every single contribution to $\cE_k^{(n)}$ (except for $\smash{\mu_k^{(n)}}$) includes a partial derivative $\partial_{i_1}\ldots \partial_{i_\ell}R_k(\rho;\lambdaT)$ for some $i_1,\ldots,i_\ell$. Expanding this in terms of $\alpha$ and $\alpha_k$ yields ratios with numerators made of products of derivatives of $\alpha$ and $\alpha_k$ and denominators amounting to some power of $\alpha$. It follows from Propositions~2.1 and 2.4 that $\cG_k^{(n)}$ can be represented in full generality by
\begin{equation}
\label{eq:Gn}
\left.\begin{array}{ll}
  \displaystyle{\cG^{(n)}_k =\, \mu_k^{(n)} +\ \blambda^{-n+1}\stackrel[\bk\le n+1]{}{\sum_{k_1=0}^{n+1}\sum_{k_2=0}^{k_1}\ldots \sum_{k_{n+1}=0}^{k_n}}
c^{(n)}_{k;k_1\ldots k_{n+1}}\frac{F_{v+2k_1}\cdot\ldots\cdot
F_{v+2k_{n+1}}}{F_v^{n+1}}}\,,\\[7.0ex]
 \bk \equiv \sum_{j=1}^{n+1} k_j\,,\end{array}\right.
\end{equation}
with the coefficient $\blambda^{-n+1}$ having been factored out for later convenience. The subscript prescription $K\le n+1$ has to be understood as a restricting condition on the possible values taken by $k_1,\ldots,k_{n+1}$. For the sake of conciseness, we refer collectively to the ratios $F_{v+2k_1}\cdot\ldots\cdot F_{v+2k_{n+1}}/F_v^{n+1}$ as the \emph{$\chi^2$-ratios} and to the coefficients $c^{(n)}_{k;k_1\ldots k_{n+1}}$ as the \emph{$c$-coefs}. Evidently the \rhs of eq.~(\ref{eq:Gn}) becomes increasingly populated as we  increase $n$. The lowest order coefficients are 
\begin{align}
\cG^{(1)}_k \ & = \ \mu_k^{(1)}\,+\,c^{(1)}_{k;00}\,+\,c^{(1)}_{k;10}\frac{F_{v+2}}{F_v}\,+\,c^{(1)}_{k;11}\frac{F_{v+2}^2}{F_v^2}\,+\,c^{(1)}_{k;20}\frac{F_{v+4}}{F_v}\,, \\[2.0ex]
\cG^{(2)}_k \ & = \ \mu_k^{(2)}\,+\,\blambda^{-1}\left( c^{(2)}_{k;000}+c^{(2)}_{k;100}\frac{F_{v+2}}{F_v} +
c^{(2)}_{k;110}\frac{F_{v+2}^2}{F_v^2}+c^{(2)}_{k;111}\frac{F_{v+2}^3}{F_v^3}\right.\nonumber\\[0.0ex]
\ & \hskip 2.3cm\left.+\ c^{(2)}_{k;200}\frac{F_{v+4}}{F_v}+c^{(2)}_{k;210}\frac{F_{v+4}F_{v+2}}{F_v^2}+c^{(2)}_{k;300}\frac{F_{v+6}}{F_v}\right)\,,\\[0.0ex]
\, & \, \, \, \vdots\nonumber
\end{align}
Without conditioning the sum to $K\le n+1$, the number of summands in eq.~(\ref{eq:Gn}) would be $2n+2 \choose n+1$ (see eq.~(1) of ref.~\cite{butler}). Owing to the restricting condition the actual number of summands is much lower. The intricacy of eq.~(\ref{eq:Gn}) is only apparent: adding separately the indices of all $\chi^2$-\cdf's at numerator and denominator and then subtracting the resulting numbers just yields $2K$. Therefore, eq.~(\ref{eq:Gn}) is just a formal way of representing a linear combination of $\chi^2$-ratios, where each $\chi^2$-\cdf has at least $v$ degrees of freedom and the overall algebraic sum of degrees of freedom amounts to $2K = 0,2,\ldots,2(n+1)$ (with denominators contributing negatively). We remark that this analytic structure is a direct consequence of the chain rule of differentiation together with the results of Propositions~2.1 and 2.4.    

Regarding the $c$-coefs, we observe that they do not depend on $\rho$ and can be only determined by explicit calculation. In spite of this, their dependence upon the Taylor coefficients of $\lambda(\epsilon)$ displays a well defined analytic structure. In order to show this, we rely upon the notions of physical and perturbative dimensions. 

\noindent {\bf{i})} We assume that $\lambda$ has physical dimension of length~$(L)$. To express this we adopt the notation $[\lambda]_\rL = L^1$. If we also assume $[\epsilon]_\rL=L^0$, then it follows $[\blambda]_\rL = [\lambda^{(1)}]_\rL = \ldots = [\lambda^{(n)}]_\rL = L^1$. Similarly, we have $[\mu]_\rL = L^1$ and $[R_k]_\rL = L^0$. Since $[\cJ]_\rL = L^0$, from eq.~(\ref{eq:finalpt}) it follows $[\cG^{(n)}]_\rL = L^1$.  Hence, eq.~(\ref{eq:Gn}) makes us conclude that
\begin{equation}
[c^{(n)}_{k;k_1\ldots k_{n+1}}]_\rL = L^{n}\,.
\label{eq:physdim}
\end{equation}
As previously explained, the $c$-coefs depend only polynomially upon the Taylor coefficients of~$\lambda(\epsilon)$. Eq.~(\ref{eq:physdim}) suggests that these polynomials are linear combinations of monomials in $\blambda$ and the components of $\lambda^{(1)},\ldots,\lambda^{(n-1)}$, each monomial having precisely degree $n$.   

\noindent {\bf{ii})} We define the perturbative dimension of a single monomial as the sum of the perturbative orders of its factors. More precisely, we let $[\blambda]_\rPT = P^0$, $[\lambda^{(1)}]_\rPT = P^1$, \ldots, $[\lambda^{(n)}]_\rPT = P^n$. Thus, for instance, we have $[\blambda^2(\lambda^{(2)}_i)^3(\lambda^{(1)}_j)^{2}]_\rPT = P^8$. We remark that $\blambda$ bears no perturbative dimension, yet it increases the physical dimension of the monomials. Since $\cG^{(n)}$ is the result of the expansion at $\Oen$, it follows that
\begin{equation}
[c^{(n)}_{k;k_1\ldots k_{n+1}}]_\rPT = P^n\,.
\label{eq:pertdim}
\end{equation}
Each monomial contributing to $c^{(n)}_{k;k_1\ldots k_{n+1}}$ has the same perturbative dimension $P^n$.

\noindent {\bf{iii})} Several monomials contributing to a given $c$-coef have the same perturbative structure and numerical prefactor and differ only by directional indices, {\it e.g.} the monomials $3\blambda^2(\lambda^{(2)}_1)^3(\lambda^{(1)}_2)^{2}$ and $3\blambda^2(\lambda^{(2)}_4)^3(\lambda^{(1)}_5)^{2}$. This is a consequence of the index structure of $\Delta_{i_1\ldots i_\ell}$: the products of Kronecker symbols contributing to the \rhs of eq.~(\ref{eq:DeltaKron}) contract the indices of the  Taylor coefficients of $\lambda(\epsilon)$ in all possible ways, thus generating an increasing number of new aggregate structures at each order of the expansion. For instance, monomials within the $c$-coefs belonging to the lowest perturbative orders can be grouped according to   
\begin{align}
\label{eq:psone}
& \Oetwo\,:\quad \zeta_1  = \sum_{i=1}^v\lambda^{(1)}_i\,,\quad \zeta_{1:2} = \sum_{i=1}^v
(\lambda^{(1)}_i)^2\,;\\[0.0ex]
\label{eq:pstwo}
& \Oethree\,:\quad \zeta_{2} = \sum_{i=1}^v\lambda^{(2)}_i\,,\quad \zeta_{12} = \sum_{i=1}^v
\lambda^{(1)}_i\lambda^{(2)}_i\,,\quad \zeta_{1:3} = \sum_{i=1}^v(\lambda^{(1)}_i)^3\,;\\[0.0ex]
& \Oefour\,:\quad \zeta_{3} = \sum_{i=1}^v\lambda^{(3)}_i\,,\quad \zeta_{13}
= \sum_{i=1}^v \lambda^{(1)}_i\lambda^{(3)}_i\,,\quad
\zeta_{2:2}=\sum_{i=1}^v(\lambda^{(2)}_i)^2\,,\nonumber\\[0.0ex]
\label{eq:psthree}
& \hskip 1.69cm
\zeta_{1:2\,2}=\sum_{i=1}^v(\lambda_i^{(1)})^2\lambda_i^{(2)}\,,\quad
\zeta_{1:4}=\sum_{i=1}^v(\lambda^{(1)}_i)^4\,;\\[0.0ex] 
&\hskip 1.1cm \vdots\nonumber
\end{align}

In view of the above considerations, we conclude that all $c$-coefs at $\Oen$ with $n\ge 2$ can be represented in full generality as linear combinations of all possible products of perturbative structures under the constraints imposed by eqs.~(\ref{eq:physdim}) and (\ref{eq:pertdim}), {\it i.e.}
\begin{equation}
c^{(n)}_{k;k_1\ldots k_{n+1}} = \sum_m \gamma^{(n,m)}_{k_1\ldots k_{n+1}}\cO^{(n)}_{k;m}\,,
\label{eq:cvsgamma}
\end{equation}
with numerical prefactors $\gamma^{(n,m)}_{k_1\ldots k_{n+1}}$ and perturbative structures $\cO^{(n)}_{k;m}$ fulfilling $[\cO^{(n)}_{k;m}]_\rL = L^n$ and $[\cO^{(n)}_{k;m}]_\rPT = P^n$. For instance, we have 
\begin{align}
\label{eq:Otwobasis}
& \cO^{(2)}_{k;m} \ \in\  \{(\lambda^{(1)}_k)^2,\,\lambda_k^{(1)}\zeta_1,\,\zeta_1^2,\,\zeta_{1:2}\}\,,
\end{align}
\begin{align}
\label{eq:Othreebasis}
& \cO^{(3)}_{k;m} \ \in\
\{(\lambda^{(1)}_k)^3,(\lambda_k^{(1)})^2\zeta_1,\,\lambda_k^{(1)}\zeta_1^2,\,\zeta_1^3,\,\lambda_k^{(1)}\zeta_{1:2},\,\zeta_{1:3},\nonumber\\[3.0ex]
& \hskip 1.69cm
\zeta_1\zeta_{1:2},\,\blambda\lambda_k^{(1)}\lambda_k^{(2)},\,\blambda\lambda_k^{(1)}\zeta_2,\,\blambda\zeta_{12},\,\blambda\lambda_k^{(2)}\zeta_1,\,\blambda\zeta_1\zeta_2\}\,,\\[2.0ex]
& \hskip 1.0cm \vdots\nonumber
\end{align}
Before working out the expansion at a given order, one should list all possible perturbative structures pertaining to that order, such as eqs.~(\ref{eq:Otwobasis}) and (\ref{eq:Othreebasis}) for $n=2,3$ respectively. A preliminary identification of all suitable structures is indeed particularly useful in order to identify groups of terms when high order calculations are performed by means of a computer algebra system (CAS), as we shall see in sect. 3.3. 

When calculating $c^{(n)}_{k;k_1\ldots k_{n+1}}$, many of the coefficients $\gamma^{(n,m)}_{k_1\ldots   k_{n+1}}$ are found to be zero. The non-vanishing ones fulfill the following property:

\begin{prop}
For $n\ge 2$ the coefficients $\gamma^{(n,m)}_{k_1\ldots k_{n+1}}$ fulfill
\begin{equation}
\stackrel[\bk\le n+1]{}{\sum_{k_1=0}^{n+1}\sum_{k_2=0}^{k_1}\ldots
  \sum_{k_{n+1}=0}^{k_n}} \gamma^{(n,m)}_{k_1\ldots k_{n+1}} = 0\,.
\label{eq:gammaconstr}
\end{equation}
\end{prop}
\begin{proof}
We first note that if $\mu\in \cD(\tau_{\rho}^{-1})$, then $\mu\in \cD(\tau_{\rho'}^{-1})$ \  $\forall\, \rho'>\rho$. Therefore, it makes sense to consider eq.~(\ref{eq:finalpt}) as $\rho\to\infty$ with $\mu$ kept fixed. In particular, we have shown previously that $\lim_{\rho\to\infty}\cJ_{kj} =~\delta_{kj}$. Moreover,  as $\rho\to\infty$ all the $\chi^2$-ratios tend to one, thus eq.~(\ref{eq:finalpt}) reduces to
\begin{align}
\lambda^{(n)}_k & = \mu_k^{(n)} + \stackrel[\bk\le n+1]{}{\sum_{k_1=0}^{n+1}\sum_{k_2=0}^{k_1}\ldots
  \sum_{k_{n+1}=0}^{k_n}} c^{(n)}_{k;k_1\ldots k_{n+1}} = \mu_k^{(n)} + \stackrel[\bk\le n+1]{}{\sum_{k_1=0}^{n+1}\sum_{k_2=0}^{k_1}\ldots
  \sum_{k_{n+1}=0}^{k_n}}\sum_{m}\gamma^{(n,m)}_{k_1\ldots
  k_{n+1}}\cO^{(n)}_{k;m} \nonumber\\[2.0ex]
& = \mu_k^{(n)}+ \sum_m \cO^{(n)}_{k,m}(\blambda,\lambda^{(1)},\ldots,\lambda^{(n-1)})\left[\stackrel[\bk\le n+1]{}{\sum_{k_1=0}^{n+1}\sum_{k_2=0}^{k_1}\ldots
  \sum_{k_{n+1}=0}^{k_n}}\gamma^{(n,m)}_{k_1\ldots k_{n+1}}\right]
\label{eq:gammazero}
\end{align}
However, in the same limit $\lambda_k \to \mu_k$, which entails order by order $\lambda_k^{(n)}\to\mu_k^{(n)}$. Hence, we infer that the sum of perturbative structures on the \rhs of eq.~(\ref{eq:gammazero}) vanishes as $\rho\to\infty$. Since in general $\cO^{(n,m)}_{k,m}(\bmu,\mu^{(1)},\ldots,\mu^{(n-1)}) \ne~0$, we conclude that eq.~(\ref{eq:gammaconstr}) is correct.
\end{proof}

The first few orders of the perturbative expansion can be worked out with little algebraic effort. Doing the calculations is useful to familiarize with the general structure discussed so far. The order $\cO(\epsilon^0)$ of the expansion has been discussed in sect.~1. We can focus on the perturbative corrections to it. From now on we assume that $\mu(\epsilon)$ has Taylor coefficients given by eq.~(\ref{eq:splittingone}). 

\subsection{Perturbative expansion at $\Oeone$}

Equations $\{\cE_k^{(1)}=0\}_{k=1}^v$ have the explicit form 
\begin{equation}
\delta\mu_k = \lambda_k^{(1)}\bR + \blambda R_k^{(1)}\,.
\end{equation}
Since $\tilde R = F_{v+2}/F_v$, the only term we need to calculate is 
\begin{equation}
R_k^{(1)} = \frac{\rd R_k}{\rd\epsilon}\biggr|_{\epsilon=0} = \sum_{j=1}^v\lambda_j^{(1)}\partial_jR_k(\rho;\lambdaT)\,.
\end{equation}
Actually, we have calculated the first order partial derivatives of $R_k$ in eq.~(\ref{eq:firstRder}). Thus, we have
\begin{equation}
\sum_{j=1}^v \cJ_{kj}\,\lambda_j^{(1)} = \delta\mu_k \,.
\label{eq:PTorderone}
\end{equation}
whence we infer $\cG^{(1)}_k = \delta\mu_k$. Choosing $\bmu=\bar\mu$ yields an important simplification:  
\begin{prop}
If $\bmu = \bar\mu$, then $\zeta_1 = 0$.
\end{prop}
\begin{proof}
It is sufficient to add side by side all eqs.~(\ref{eq:PTorderone}) for $k=1,\ldots,v$ to get 
\begin{equation}
\frac{\zeta_1}{2}\left[(v+2)\frac{F_{v+4}}{F_{v}}-v\frac{F^2_{v+2}}{F_v^2}\right] =
\sum_{k=1}^v\delta\mu_k = 0\,.
\end{equation}
Since the quantity in square brackets is strictly positive, we conclude that $\zeta_1=0$.   
\end{proof}
As can be easily understood, having $\zeta_1=0$ results in a huge simplification of the algebra. Indeed, $\zeta_1$ belongs to many perturbative structures contributing to $\smash{\cG^{(n)}}$ for $n\ge 2$. For instance, the set of structures given in eq.~(\ref{eq:Otwobasis}) is reduced to only two elements in place of four when $\zeta_1=0$, while the one given in eq.~(\ref{eq:Othreebasis}) is reduced to six elements in place of twelve. 

\subsection{Perturbative expansion at $\Oetwo$}

The subleading correction $\lambda^{(2)}$ is obtained from the equations
\begin{equation}
0 = \lambda_k^{(2)}\bR + \lambda_k^{(1)}R_k^{(1)} + \blambda R_k^{(2)}\,.
\end{equation}
Most of the contributions to the three terms on the \rhs are worked out easily at this point. For instance, we have
\begin{equation}
\lambda_k^{(1)}R_k^{(1)} =
\frac{1}{\blambda}\left[2(\lambda_k^{(1)})+\lambda_k^{(1)}\zeta_1\right]\frac{F_{v+4}}{F_v}
- \frac{\lambda_k^{(1)}\zeta_1}{\blambda}\frac{F_{v+2}^2}{F_v^2}- \frac{(\lambda_k^{(1)})^2}{\blambda}\frac{F_{v+2}}{F_v}\,,
\end{equation}
\begin{equation}
\lambda_k^{(2)}\bR+\blambda R_k^{(2)} = \frac{1}{2}\frac{\rd^2 R_k}{\rd \epsilon^2}\biggr|_{\epsilon=0} =
\sum_{j=1}^v\cJ_{kj}\lambda^{(2)}_j + \frac{\blambda}{2}\sum_{j_1j_2=1}^v\lambda^{(1)}_{j_1}\lambda^{(1)}_{j_2}\partial_{j_1j_2}R_k(\rho;\lambdaT)\,.
\end{equation}
To keep things general, we make no assumptions on $\bmu$ here. Accordingly, we retain all $\chi^2$-ratios proportional to powers of $\zeta_1$. The evaluation of the second derivatives $\partial_{j_1j_2}R_k(\rho;\lambdaT)$ requires a few pages of tedious algebraic work, which we cannot detail. The result of the calculations is given by
\begin{align}
\sum_{j=1}^v \cJ_{kj}\lambda_j^{(2)} & = - \frac{1}{8\blambda}\left[\zeta_1^2
  + 2\zeta_{1:2} + 4\lambda_k^{(1)}\zeta_1 +
  8(\lambda_k^{(1)})^2\right]\frac{F_{v+6}}{F_v}\nonumber\\[0.8ex]
& + \frac{1}{8\blambda}\left[3\zeta_1^2 + 2\zeta_{1:2}
+ 4\lambda_k^{(1)}\zeta_1\right]\frac{F_{v+4}F_{v+2}}{F_v^2}
-
\frac{\zeta_1^2}{4\blambda}\frac{F_{v+2}^3}{F_v^3}\nonumber\\[2.4ex]
& +\frac{1}{2\blambda}[\zeta_{1:2} + 2(\lambda_k^{(1)})^2]\frac{F_{v+4}}{F_v} -\frac{\zeta_{1:2}}{2\blambda}\frac{F_{v+2}^2}{F_v^2}\,.
\label{eq:Gntwo}
\end{align}
\begin{table}[!t]
  \small
  \begin{center}
    \begin{tabular}{l|cccccc}
      \hline\hline\\[-2.0ex]
      & $\dfrac{F_{v+6}}{F_v}$ & $\dfrac{F_{v+4}F_{v+2}}{F_v^2}$ &
      $\dfrac{F_{v+2}^3}{F_v^3}$ &
      $\dfrac{F_{v+4}}{F_v}$ & $\dfrac{F_{v+2}^2}{F_{v}^2}$ \\ \\[-2.5ex]
      \hline\\[-2.5ex]
      $(\loko)^2$      &  $-1$ &  $\phantom{-}0$ & $\phantom{-}0$
      &  $\phantom{-}1$ & $\phantom{-}0$ \\[1.0ex]
      $\loko\zeta_1$ & $-\dfrac{1}{2}$ & $\phantom{-}\dfrac{1}{2}$ &
      $\phantom{-}0$ & $\phantom{-}0$ & $\phantom{-}0$ \\[2.0ex]
      $\zeta_1^2$ & $-\dfrac{1}{8}$ & $\phantom{-}\dfrac{3}{8}$ &
      $-\dfrac{1}{4}$ & $\phantom{-}0$ & $\phantom{-}0$ \\[2.0ex]
      $\zeta_{1:2}$  &  $-\dfrac{1}{4}$  &  $\phantom{-}\dfrac{1}{4}$  &
      $\phantom{-}0$ & $\phantom{-}\dfrac{1}{2}$
      & $-\dfrac{1}{2}$ \\[2.0ex]
      \hline\hline
    \end{tabular}
    \vskip 0.1cm
    \caption{Coefficients $\gamma^{(2,m)}_{k_1k_2 k_3}$.\label{tab:coeftwo}} 
  \end{center}
  \vskip -0.5cm
\end{table}

\noindent We notice that all four structures listed in eq.~(\ref{eq:Otwobasis}) contribute to the $c$-coefs $c^{(2)}_{k;k_1k_2k_3}$. 

In Table~\ref{tab:coeftwo} we collect the coefficients $\gamma^{(2,m)}_{k_1k_2k_3}$. Instead of naming rows and columns respectively according to the values of $m$ and the triples $(k_1,k_2,k_3)$, for the sake of readability we identify each table entry by the perturbative-structure and the $\chi^2$-ratio it refers to. This way of tabulating coefficients becomes particularly informative at higher orders. We observe that adding the entries of each table row yields zero, in accordance with eq.~(\ref{eq:gammaconstr}). 

\subsection{Perturbative expansion at higher orders}

Paper-and-pencil calculations become prohibitively expensive at higher orders. Fortunately, it is not difficult to work out the algebra with the assistance of a CAS. For the reader's convenience, in App.~C we attach some essential and correctly working \Maple procedures, which help work out the algebra. The code is split into three blocks, that we shortly review. 

The first code block (C.1) contains a procedure {\tt Delta()}, which computes the coefficient $\Delta_{i_1\ldots i_n}$. The procedure argument is assumed to be a list of {\tt nonnegint} items; alternatively the procedure returns unevaluated. The input list is first sorted in ascending order, then the  multiplicity set is identified. The procedure computes the \rhs of eq.~(\ref{eq:deltasymbol}) and returns its numerical value. 

The second code block (C.2) performs the algebraic work related to the perturbative expansion of eqs.~(\ref{eq:truncspectrum}). Before submitting it to evaluation, the user is assumed to assign a {\tt nonnegint} variable~{\tt v} representing the number of dimensions, and a {\tt nonnegint} variable {\tt n $\le$ 4} representing the highest perturbative order processed by the program.  The code block starts with a pair of procedures, {\tt DerAlpha()} and {\tt   DerAlphak()}, which encode respectively  eqs.~(\ref{eq:deralpha}) and (\ref{eq:deralphak}). Then, it performs a Taylor expansion of the \rhs of eqs.~(\ref{eq:truncspectrum}) up to ${\rm O}(\epsilon^{\tt n})$. Taylor coefficients are stored within {\tt indexable} objects {\tt h[j,k]}, the indices {\tt j} and {\tt k} representing respectively the perturbative order and the physical direction. At this stage, {\tt h[j,k]} includes a sum of potentially many terms. The summands contain derivatives of $\alpha$ and $\alpha_k$, which are purely symbolic objects at this stage. Their evaluation requires sequences of prescriptions, stored within the variables {\tt C0A}, {\tt C0Ak},\ldots, {\tt C4A}, {\tt C4Ak}. Algebraic simplifications are performed in the last few lines, where partial results are stored within {\tt indexable} objects {\tt h00},\ldots, {\tt h10}, so as to allow for an offline analysis of the single steps.

\begin{table}[!t]
  \small
  \begin{center}
    \begin{tabular}{l|ccccccc}
      \hline\hline\\[-2.0ex]
      & $\dfrac{F_{v+8}}{F_v}$ & $\dfrac{F_{v+6}F_{v+2}}{F_v^2}$ & $\dfrac{F_{v+4}^2}{F_v^2}$ &
      $\dfrac{F_{v+6}}{F_{v}}$ & $\dfrac{F_{v+4}F_{v+2}}{F_v^2}$ &
      $\dfrac{F_{v+4}}{F_v}$ & $\dfrac{F_{v+2}^2}{F_v^2}$ \\ \\[-2.5ex]
      \hline\\[-2.0ex]
      $(\loko)^3$                                    & $-1$   & $\phantom{-}0$     & $\phantom{-}0$    & $\phantom{-}2$   & $\phantom{-}0$    & $-1$   & $\phantom{-}0$  \\[1.0ex]
      $\loko\zeta_{1:2}$                           & $-\dfrac{1}{4}$ & $\phantom{-}0$    & $\phantom{-}\dfrac{1}{4}$ & $\phantom{-}\dfrac{1}{2}$ & $-\dfrac{1}{2}$ & $\phantom{-}0$   & $\phantom{-}0$  \\[2.0ex]
      $\zeta_{1:3}$                                & $-\dfrac{1}{6}$ & $\phantom{-}\dfrac{1}{6}$ & $\phantom{-}0$    & $\phantom{-}\dfrac{1}{2}$ & $-\dfrac{1}{2}$ & $-\dfrac{1}{2}$ & $\phantom{-}\dfrac{1}{2}$ \\[2.5ex]
      $\blambda\loko\ltko$                        & $\phantom{-}0$   & $\phantom{-}0$     & $\phantom{-}0$    & $-2$    & $\phantom{-}0$    & $\phantom{-}2$  & $\phantom{-}0$  \\[1.5ex]
      $\blambda\loko\zeta_2$                & $\phantom{-}0$   & $\phantom{-}0$     & $\phantom{-}0$    & $-\dfrac{1}{2}$  & $\phantom{-}\dfrac{1}{2}$ & $\phantom{-}0$   & $\phantom{-}0$  \\[2.0ex]
      $\blambda\zeta_{12}$  & $\phantom{-}0$   & $\phantom{-}0$     & $\phantom{-}0$    & $-\dfrac{1}{2}$  & $\phantom{-}\dfrac{1}{2}$ & $\phantom{-}1$  & $-1$  \\[1.5ex]
      \hline\hline
    \end{tabular}
    \vskip 0.0cm
    \caption{Coefficients $\gamma^{(3,m)}_{k_1k_2 k_3 k_4}$. We assume in this
      case $\zeta_1=0$.\label{tab:coefthree}}
  \end{center}
\end{table}

The third code block (C.3) illustrates in a specific case a numerical technique which we have devised for the determination of the coefficients $\gamma^{(n,m)}_{k_1\ldots k_{n+1}}$. The code processes the coefficients corresponding to $n=3$ and $(k_1,k_2,k_3,k_4) = (3,0,0,0)$, {\it i.e.} those entering the $c$-coef multiplying the $\chi^2$-ratio $F_{v+6}/F_v$. It also assumes $\bmu = \bar\mu$. As previously discussed, this simplifies the basis of perturbative structures to
\begin{equation}
\cO^{(3)}_{k;m} \ \in\
\{(\lambda^{(1)}_k)^3,\,\lambda_k^{(1)}\zeta_{1:2},\,\zeta_{1:3},\,
\blambda\lambda_k^{(1)}\lambda_k^{(2)},\,\blambda\lambda_k^{(1)}\zeta_2,\,\blambda\zeta_{12}\}\,.
\end{equation}
The algebraic sum pointed to by {\tt h10[3,k]} at the end of the second code block has no knowledge of these structures. In order to identify them within {\tt h10[3,k]}, we need to group terms properly. Instead of proceeding at an algebraic level, which would be computationally demanding, we adopt a numerical approach, based on the use of eq.~(\ref{eq:cvsgamma}) as a square linear system fulfilled by the coefficients $\smash{\gamma^{(n,m)}_{k_1\ldots k_{n+1}}}$. Having subtracted from {\tt h10[3,k]} all contributions appearing in  eq.~(\ref{eq:termsln}), we extract from it all terms proportional to $F_{v+6}/F_v$, whose sum amounts to $\smash{-c^{(3)}_{k;3000}}$. Then, for each $k$ we  assign $\lambda^{(1)}$ and $\lambda^{(2)}$ random values (chosen so that $\zeta_1=0$), from which we compute $\smash{\cO^{(3)}_{k;1}}$, \ldots, $\smash{\cO^{(3)}_{k;6}}$ and $\smash{c^{(3)}_{k;3000}}$. The random matrix $\smash{\cO^{(3)}_{k;m}}$ thus obtained is non-singular, hence eq.~(\ref{eq:termsln}) can be solved with respect to $\smash{\gamma^{(3,m)}_{3000}}$.  The solution is independent of the random numbers generated. This represents a strong signal that our determination is correct, yet a real check consists of an algebraic comparison between the reconstructed coefficient $c^{(3)}_{k;3000}$ and the one extracted from {\tt h10[3,k]}. In Tables~\ref{tab:coefthree} and \ref{tab:coeffour} we report the coefficients $\smash{\gamma^{(3,m)}_{k_1k_2k_3k_4}}$ and $\smash{\gamma^{(4,m)}_{k_1k_2k_3k_4k_5}}$ under the assumption $\zeta_1=0$.     

\section{Properties of the first few perturbative coefficients}

So far we have focused on formal aspects of the perturbative expansion of $\tau_\rho^{-1}$ with the aim of proving its theoretical and computational feasibility. To establish the level of accuracy reached in approximating the reconstruction operator by a handful of perturbative contributions, we need to investigate some analytic properties of the perturbative coefficients of~$\lambda$.

Hereinafter we consider perturbative series truncated at different orders, for which we adopt the notation
\begin{equation}
\lambdabar^{(n)}_k \,\equiv\, \sum_{i=0}^n \lambda_k^{(i)}\,.
\end{equation}
Fig.~\ref{fig:evaex} shows an illustrative example of perturbative reconstructions at $v=4$. We concentrate on it in the present and next sections. The plots in Fig.~\ref{fig:evaex} have been produced as follows. First of all, in order to test perturbation theory on eigenvalues characterized by a relatively large ratio $\lambda_\text{max}/\lambda_\text{min}$, we have chosen $\lambdaEX = \{0.1,0.3,0.8,2.2\}$ as the eigenvalue set to reconstruct (accordingly, we have $\lambda_\text{max}/\lambda_\text{min} = 22\gg 1$, yielding a highly asymmetric Gaussian ellipsoid). By means of numerical techniques detailed in ref.~\cite{palombi4}, we have then computed $\mu = \tau_\rho\cdot\lambdaEX$ for several values of $\rho$. In correspondence with each pair $(\rho,\mu)$ we have finally reconstructed $\lambda$ up to the fourth perturbative order, having chosen in all cases $\bmu = \bar\mu$. We notice from the plots that in most cases the error made by truncating the expansion at a given order increases for lower values of $\rho$. Moreover, the error is larger for eigenvalues at the extremes of the eigenvalue set and milder in the center of it. We also notice that the convergence pattern for the lower eigenvalues is radically different than for the higher ones. Indeed, the perturbative series of $\lambda_1$ and $\lambda_2$ converges with alternate signs, while the perturbative series of $\lambda_3$ and $\lambda_4$ displays a monotonic character. 

In order to explain the observed behavior, we first concentrate on the leading contribution~$\blambda$. If $(\rho^*,\bmu)$ fulfills the constraint $0<\bmu< \rho^*/(v+2)$, a solution $\blambda$ to eq.~(\ref{eq:tspectrum}) exists for all pairs $(\rho,\bmu)$ with $\rho>\rho^*$. We can consider $\blambda$ as a function of $\rho$ at fixed $\bmu$. The analytic form of $\cT_\rho$ tells us that $\lim_{\rho\to\infty}\blambda = \bmu$. Moreover, we know that $\bar\mu< \rho^*/(v+2)$ provided $\mu\in\cD(\tau_{\rho^*}^{-1})$. Since $\mu_1\le\bar\mu\le\mu_v$, by continuity we conclude that  
\begin{equation}
\exists\, \hat\rho\,:\qquad \cT_\rho^{-1}(\mu_1)\le \blambda \le \cT_\rho^{-1}(\mu_{v}) \qquad \forall\ \rho\ge\hat\rho\,.
\label{eq:blambdalocone}
\end{equation}
At sufficiently large $\rho$ the above inequality holds true with $\mu_1$ and $\mu_v$ being respectively replaced by $\mu_i$ and $\mu_{i+1}$, where $i$ is such that $\mu_i\le\bar\mu\le\mu_{i+1}$. 

\begin{center}
  \begin{landscape}
    \begin{table}
      \vskip -0.5cm
      \small
      \addtolength{\tabcolsep}{-2pt}
      \begin{tabular}{l|ccccccccccccc}
        \hline\hline\\[0.0ex]
        & $\dfrac{F_{v+10}}{F_v}$
        & $\dfrac{F_{v+8}F_{v+2}}{F_v^2}$
        & $\dfrac{F_{v+6}F_{v+4}}{F_v^2}$
        & $\dfrac{F_{v+4}^2F_{v+2}}{F_v^3}$
        & $\dfrac{F_{v+8}}{F_v}$
        & $\dfrac{F_{v+6}F_{v+2}}{F_v^2}$
        & $\dfrac{F_{v+4}^2}{F_v^2}$
        & $\dfrac{F_{v+4}F_{v+2}^2}{F_v^3}$
        & $\dfrac{F_{v+6}}{F_v}$
        & $\dfrac{F_{{v+4}}F_{{v+2}}}{{F_{{v}}}^{2}}$
        & $\dfrac{F_{v+2}^3}{F_{v}^3}$
        & $\dfrac{F_{v+4}}{F_v}$
        & $\dfrac{F_{v+2}^2}{F_v^2}$\\[2.0ex]
        \hline\\[0.0ex]
$(\lambda^{(1)}_k)^4$ & $-1$ & $\phantom{-}0$ & $\phantom{-}0$ & $\phantom{-}0$ & $\phantom{-}3$ & $\phantom{-}0$ & $\phantom{-}0$ & $\phantom{-}0$ & $-3$ & $\phantom{-}0$ & $\phantom{-}0$ & $\phantom{-}1$ & $\phantom{-}0$  \\[2.0ex]
$(\lambda^{(1)}_k)^2\zeta_{1:2}$ & $-1/4$ & $\phantom{-}0$ & $\phantom{-}1/4$ & $\phantom{-}0$ & $\phantom{-}3/4$ & $-1/2$ & $-1/4$ & $\phantom{-}0$ & $-1/2$ & $\phantom{-}1/2$ & $\phantom{-}0$ & $\phantom{-}0$ & $\phantom{-}0$  \\[2.0ex]
$\lambda^{(1)}_k\zeta_{1:3}$ & $-1/6$ & $\phantom{-}0$ & $\phantom{-}1/6$ & $\phantom{-}0$ & $\phantom{-}1/2$ & $\phantom{-}0$ & $-1/2$ & $\phantom{-}0$ & $-1/2$ & $\phantom{-}1/2$ & $\phantom{-}0$ & $\phantom{-}0$ & $\phantom{-}0$  \\[2.0ex]
$\zeta_{1:2}^2$ & $-1/32$ & $\phantom{-}1/32$ & $\phantom{-}1/16$ & $-1/16$ & $\phantom{-}1/8$ & $-1/4$ & $-1/8$ & $\phantom{-}1/4$ & $-1/8$ & $\phantom{-}3/8$ & $-1/4$ & $\phantom{-}0$ & $\phantom{-}0$  \\[2.0ex]
$\zeta_{1:4}$ & $-1/8$ & $\phantom{-}1/8$ & $\phantom{-}0$ & $\phantom{-}0$ & $\phantom{-}1/2$ & $-1/2$ & $\phantom{-}0$ & $\phantom{-}0$ & $-3/4$ & $\phantom{-}3/4$ & $\phantom{-}0$ & $\phantom{-}1/2$ & $-1/2$  \\[2.0ex]
$\blambda(\lambda^{(1)}_k)^2\lambda^{(2)}_k$ & $\phantom{-}0$ & $\phantom{-}0$ & $\phantom{-}0$ & $\phantom{-}0$ & $-3$ & $\phantom{-}0$ & $\phantom{-}0$ & $\phantom{-}0$ & $\phantom{-}6$ & $\phantom{-}0$ & $\phantom{-}0$ & $-3$ & $\phantom{-}0$  \\[2.0ex]
$\blambda(\lambda^{(1)}_k)^2\zeta_2$ & $\phantom{-}0$ & $\phantom{-}0$ & $\phantom{-}0$ & $\phantom{-}0$ & $-1/2$ & $\phantom{-}1/2$ & $\phantom{-}0$ & $\phantom{-}0$ & $\phantom{-}1/2$ & $-1/2$ & $\phantom{-}0$ & $\phantom{-}0$ & $\phantom{-}0$  \\[2.0ex]
$\blambda\lambda^{(2)}_k\zeta_{1:2}$ & $\phantom{-}0$ & $\phantom{-}0$ & $\phantom{-}0$ & $\phantom{-}0$ & $-1/4$ & $\phantom{-}0$ & $\phantom{-}1/4$ & $\phantom{-}0$ & $\phantom{-}1/2$ & $-1/2$ & $\phantom{-}0$ & $\phantom{-}0$ & $\phantom{-}0$  \\[2.0ex]
$\blambda\zeta_{1:2}\zeta_2$ & $\phantom{-}0$ & $\phantom{-}0$ & $\phantom{-}0$ & $\phantom{-}0$ & $-1/8$ & $\phantom{-}1/4$ & $\phantom{-}1/8$ & $-1/4$ & $\phantom{-}1/4$ & $-3/4$ & $\phantom{-}1/2$ & $\phantom{-}0$ & $\phantom{-}0$  \\[2.0ex]
$\blambda\lambda^{(1)}_k\zeta_{12}$ & $\phantom{-}0$ & $\phantom{-}0$ & $\phantom{-}0$ & $\phantom{-}0$ & $-1/2$ & $\phantom{-}0$ & $\phantom{-}1/2$ & $\phantom{-}0$ & $\phantom{-}1$ & $-1$ & $\phantom{-}0$ & $\phantom{-}0$ & $\phantom{-}0$  \\[2.0ex]
$\blambda\zeta_{1:2\, 2}$ & $\phantom{-}0$ & $\phantom{-}0$ & $\phantom{-}0$ & $\phantom{-}0$ & $-1/2$ & $\phantom{-}1/2$ & $\phantom{-}0$ & $\phantom{-}0$ & $\phantom{-}3/2$ & $-3/2$ & $\phantom{-}0$ & $-3/2$ & $\phantom{-}3/2$  \\[2.0ex]
${\blambda}^2(\lambda^{(2)}_k)^2$ & $\phantom{-}0$ & $\phantom{-}0$ & $\phantom{-}0$ & $\phantom{-}0$ & $\phantom{-}0$ & $\phantom{-}0$ & $\phantom{-}0$ & $\phantom{-}0$ & $-1$ & $\phantom{-}0$ & $\phantom{-}0$ & $\phantom{-}1$ & $\phantom{-}0$  \\[2.0ex]
${\blambda}^2\lambda^{(2)}_k\zeta_2$ & $\phantom{-}0$ & $\phantom{-}0$ & $\phantom{-}0$ & $\phantom{-}0$ & $\phantom{-}0$ & $\phantom{-}0$ & $\phantom{-}0$ & $\phantom{-}0$ & $-1/2$ & $\phantom{-}1/2$ & $\phantom{-}0$ & $\phantom{-}0$ & $\phantom{-}0$  \\[2.0ex]
${\blambda}^2\zeta_2^2$ & $\phantom{-}0$ & $\phantom{-}0$ & $\phantom{-}0$ & $\phantom{-}0$ & $\phantom{-}0$ & $\phantom{-}0$ & $\phantom{-}0$ & $\phantom{-}0$ & $-1/8$ & $\phantom{-}3/8$ & $-1/4$ & $\phantom{-}0$ & $\phantom{-}0$  \\[2.0ex]
${\blambda}^2\zeta_{2:2}$ & $\phantom{-}0$ & $\phantom{-}0$ & $\phantom{-}0$ & $\phantom{-}0$ & $\phantom{-}0$ & $\phantom{-}0$ & $\phantom{-}0$ & $\phantom{-}0$ & $-1/4$ & $\phantom{-}1/4$ & $\phantom{-}0$ & $\phantom{-}1/2$ & $-1/2$  \\[2.0ex]
${\blambda}^2\lambda^{(1)}_k\lambda^{(3)}_k$ & $\phantom{-}0$ & $\phantom{-}0$ & $\phantom{-}0$ & $\phantom{-}0$ & $\phantom{-}0$ & $\phantom{-}0$ & $\phantom{-}0$ & $\phantom{-}0$ & $-2$ & $\phantom{-}0$ & $\phantom{-}0$ & $\phantom{-}2$ & $\phantom{-}0$  \\[2.0ex]
${\blambda}^2\lambda^{(1)}_k\zeta_3$ & $\phantom{-}0$ & $\phantom{-}0$ & $\phantom{-}0$ & $\phantom{-}0$ & $\phantom{-}0$ & $\phantom{-}0$ & $\phantom{-}0$ & $\phantom{-}0$ & $-1/2$ & $\phantom{-}1/2$ & $\phantom{-}0$ & $\phantom{-}0$ & $\phantom{-}0$  \\[2.0ex]
${\blambda}^2\zeta_{13}$ & $\phantom{-}0$ & $\phantom{-}0$ & $\phantom{-}0$ & $\phantom{-}0$ & $\phantom{-}0$ & $\phantom{-}0$ & $\phantom{-}0$ & $\phantom{-}0$ & $-1/2$ & $\phantom{-}1/2$ & $\phantom{-}0$ & $\phantom{-}1$ & $-1$  \\[2.0ex]
        \hline\hline
      \end{tabular}
      \vskip 0.4cm
      \caption{Coefficients $\gamma^{(4,m)}_{k_1k_2 k_3 k_4 k_5}$. We
      assume here $\zeta_1=0$.\label{tab:coeffour}\\[50.0cm]}
  \end{table}
\end{landscape}
\end{center}
Eq.~(\ref{eq:blambdalocone}) does not tell where $\blambda$ is placed in relation to the full eigenvalue spectrum. To find such an estimate, we can resort to eq.~(24) of ref.~\cite{palombi4}. Based on arguments which are completely analogous to those used therein, we arrive easily at 
\begin{equation}
\frac{\rho}{2v+1}\frac{M\bigl(v,v+3/2,\rho/(2\lambda_1)\bigr)}{M\bigl(v,v+1/2,\rho/(2\lambda_1)\bigr)}
\,\le\, \cT_\rho(\blambda)\, \le\, \frac{\rho}{3}
\frac{M\bigl(1,5/2,\rho/(2\lambda_v)\bigr)}{M\bigl(1,3/2,\rho/(2\lambda_v)\bigr)}\,,
\end{equation}
with $M(a,b,z)$ denoting the Kummer function
\begin{equation}
M(a,b,z) = \sum_{n=0}^{\infty}\frac{1}{n!}\frac{(a)_n}{(b)_n}z^n\,,\qquad
(x)_n = \frac{\Gamma(x+n)}{\Gamma(x)}\,.
\end{equation}

\begin{center}
\begin{figure}[!t]
\begin{center}
\includegraphics[width=0.89\textwidth]{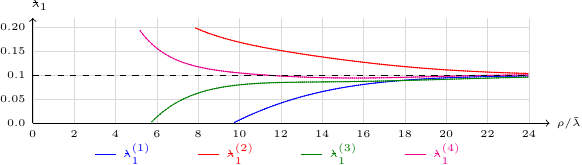}
\vskip 0.4cm
\includegraphics[width=0.89\textwidth]{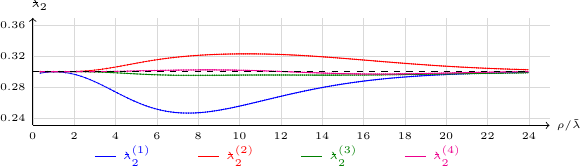}
\vskip 0.4cm
\includegraphics[width=0.89\textwidth]{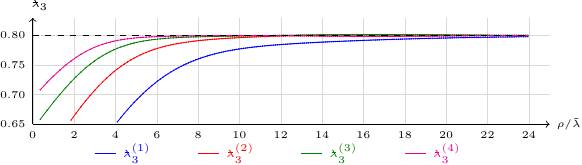}
\vskip 0.4cm
\includegraphics[width=0.89\textwidth]{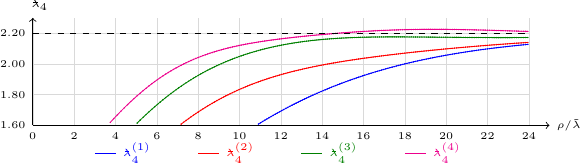}
\vskip 0.3cm
\caption{\small An example of perturbative eigenvalue reconstruction at $v=4$. The unconstrained spectrum $\lambdaEX = \{0.1,0.3,0.8,2.2\}$ is represented by black dashed lines. The solid lines correspond to the four levels of approximation obtained by truncating the perturbative series at the first to fourth order.\label{fig:evaex}}    
\end{center}
\end{figure}
\end{center}
\vskip -0.7cm
As a consequence of the asymptotic limit $\cT_\rho(\blambda)\,\widesim{\rho\to\infty}\,\blambda$ and (see {\it e.g.} chap.~13 of ref.~\cite{abramowitz}), we have
\begin{align}
& \frac{\rho}{2v+1}\frac{M\bigl(v,v+3/2,\rho/(2\lambda_1)\bigr)}{M\bigl(v,v+1/2,\rho/(2\lambda_1)\bigr)}
\ \widesim{\rho\to\infty}{} \ \lambda_1\,, \\[2.0ex]
& \frac{\rho}{3}\frac{M\bigl(1,5/2,\rho/(2\lambda_v)\bigr)}{M\bigl(1,3/2,\rho/(2\lambda_v)\bigr)} 
\ \widesim{\rho\to\infty}{} \ \lambda_v\,,
\end{align}
we conclude that if $\rho$ is sufficiently large, then $\lambda_1\le \blambda\le \lambda_v$. As intuitively expected, non--linear effects are mitigated in the region of weak truncation. The situation is qualitatively depicted in~Fig.~\ref{fig:levels}.

We now consider the first few corrections to $\blambda$. In general, as far as we are concerned with their numerical computation, we can limit ourselves to solve eqs.~(\ref{eq:Gn}) one after another by means of a linear solver. Nevertheless, the matrix structure of $\cJ^{-1}$ is simple. The off--diagonal entries are all the same, independently of the matrix indices. This allows us to invert eqs.~(\ref{eq:Gn}) analytically. When the expression thus obtained has several contributions, we hardly find it better than its numerical computation. Yet, this is not the case with the first few perturbative corrections, of which we want to estimate the range of variation. 

\begin{center}
\begin{figure}[!t]
\begin{center}
\includegraphics[width=0.7\textwidth]{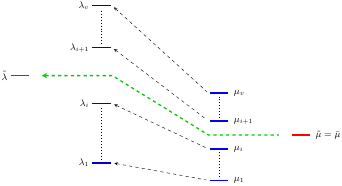}
\vskip 0.0cm
\caption{\small A schematic diagram showing the action of the operators $\tau_\rho^{-1}$ (black dashed lines) and $\cT_\rho^{-1}$ (green thick dashed line) when $\bmu=\bar\mu$ and $\mu\in\cD(\tau_{\rho}^{-1})$. If $1\le i\le v-1$ is such that  $\mu_i\le\bar\mu\le\mu_{i+1}$, then for sufficiently large $\rho$ we also have $\lambda_i\le\blambda\le\lambda_{i+1}$.\label{fig:levels}}      
\end{center}
\end{figure}
\end{center}
\vskip -0.25cm

Since $\cO^{(n)}_{k;m}$ is either index--free or dependent upon $k$ via monomials 
$(\lambda_k^{(i_1)})^{m_1}\ldots(\lambda_k^{(i_r)})^{m_r}$ with
$i_1m_1+\ldots+i_rm_r\le n$, the only algebraic ingredients we need for the
analytic inversion of eqs.~(\ref{eq:Gn}) are the sums
\begin{alignat}{3}
& \text{{\it i})} & \qquad & \sum_{k=1}^v (\cJ^{-1})_{jk} = 2\cD^{-1}\,,\qquad \cD\,\equiv\,
(v+2)\frac{F_{v+4}}{F_v} - v\frac{F_{v+2}^2}{F_v^2}\,,\\[2.0ex]
& \text{{\it ii})}& \qquad &
\sum_{j=1}^v(\cJ^{-1})_{kj}(\lambda^{(i_1)}_j)^{m_1}\ldots(\lambda^{(i_r)}_j)^{m_r}
 \nonumber\\[0.0ex]
 & & & \hskip 1.0cm =  (\lambda^{(i_1)}_k)^{m_1}\ldots(\lambda^{(i_r)}_k)^{m_r}\frac{F_v}{F_{v+4}}
- \zeta_{i_1:m_1\ldots i_r:m_r}\cD^{-1}\left(\frac{F_{v+4}}{F_v}-\frac{F_{v+2}^2}{F_v^2}\right)\,,
\end{alignat}
where $\zeta_{i_1:m_1\ldots i_r:m_r}\equiv\sum_{j=1}^v
(\lambda^{(i_1)}_j)^{m_1}\ldots(\lambda^{(i_r)}_j)^{m_r}$ is defined in analogy with eqs.~(\ref{eq:psone})--(\ref{eq:psthree}). We are now ready to work out the algebra. In first place, a straightforward calculation yields 
\begin{equation}
\lambda^{(1)}_k = \frac{F_v}{F_{v+4}}(\mu_k - \bar \mu) -2\cD^{-1}(\bmu-\bar\mu)\,.
\end{equation}
Hence, we infer that 
\begin{equation}
\sign(\lambda_k^{(1)}) = \sign(\mu_k-\bar\mu) \qquad \text{if}\ \  \bmu=\bar\mu\,.
\label{eq:firstordsign}
\end{equation}
We recognize that the first perturbative correction to the leading term $\blambda$ is positive for the higher eigenvalues and negative for the lower ones. 

The algebraic evaluation of the second perturbative correction to $\blambda$ is as easy as the first one, yet estimating its range of variation is somewhat more difficult. In this case we let $\bmu = \bar\mu$ from the very beginning. Under this assumption, we see that the \rhs of eq.~(\ref{eq:Gntwo}) reduces to a linear combination of four terms, namely
\begin{alignat}{3}
\cG^{(2)}_k & = \sum_{m=1}^4 \cG^{(2,m)}_k\,,& &\\[1.5ex]
\cG^{(2,1)}_k & = - \frac{1}{4\blambda}\left[
  \zeta_{1:2} +
  4(\lambda_k^{(1)})^2\right]\frac{F_{v+6}}{F_v}\,,\qquad & \cG^{(2,2)}_k & =
\frac{\zeta_{1:2} }{4\blambda}\frac{F_{v+4}F_{v+2}}{F_v^2},\\[1.5ex]
\cG^{(2,3)}_k & = \frac{1}{2\blambda}[\zeta_{1:2} +
2(\lambda_k^{(1)})^2]\frac{F_{v+4}}{F_v}\,,\qquad & \cG^{(2,4)} & = -\frac{\zeta_{1:2}}{2\blambda}\frac{F_{v+2}^2}{F_v^2}\,.
\end{alignat}
Correspondingly, the solution of eq.~(\ref{eq:Gntwo}) is the sum of four contributions: 
\begin{alignat}{3}
\lambda^{(2)}_k & = \sum_{m=1}^4\lambda^{(2,m)}_k\,,& &\\[1.5ex]
\lambda^{(2,1)}_k & =
-\frac{(\lambda_k^{(1)})^2}{\blambda}\frac{F_{v+6}}{F_{v+4}} +
\frac{\zeta_{1:2}}{2\blambda}\frac{F_{v+6}}{F_{v+4}}\cD^{-1}\left(\frac{F_{v+4}}{F_v}-2\frac{F_{v+2}^2}{F_v^2}\right)\,,
& &
\end{alignat}
\begin{alignat}{3}
\lambda^{(2,2)}_k & =
\frac{\zeta_{1:2}}{2\blambda}\cD^{-1}\frac{F_{v+4}F_{v+2}}{F_v^2}\,, \\[2.0ex]
\lambda^{(2,3)}_k & = \frac{(\lambda^{(1)}_k)^2}{\blambda} +
\frac{\zeta_{1:2}}{\blambda}\cD^{-1}\frac{F_{v+2}^2}{F_v^2}\,,\\[2.0ex]
\lambda^{(2,4)} & = -\frac{\zeta_{1:2}}{\blambda}\cD^{-1}\frac{F_{v+2}^2}{F_v^2}\,.
\end{alignat}
Hence, we have
\begin{equation}
\lambda^{(2)}_k =
\frac{(\lambda_k^{(1)})^2}{\blambda}\left(1-\frac{F_{v+6}}{F_{v+4}}\right) +
\frac{\zeta_{1:2}}{2\blambda}\cD^{-1}\left(\frac{F_{v+6}}{F_v}+\frac{F_{v+4}F_{v+2}}{F_v^2}
- 2\frac{F_{v+6}}{F_{v+4}}\frac{F_{v+2}^2}{F_v^2}\right)\,.
\end{equation}
Since $F_{v+6}<F_{v+4}$, the first contribution to the \rhs is certainly positive. As for the second one, we define
\begin{equation}
\Xi = \frac{F_{v+6}}{F_v}+\frac{F_{v+4}F_{v+2}}{F_v^2}
- 2\frac{F_{v+6}}{F_{v+4}}\frac{F_{v+2}^2}{F_v^2}\,.
\label{eq:xieq}
\end{equation}
Plots reported in Fig.~\ref{fig:Fineq} show that $\Xi\ge 0$. Hence, we conclude that $\lambda_k^{(2)}\ge 0$. In other words, all the eigenvalues receive a positive contribution from the second perturbative correction. Together with eq.~(\ref{eq:firstordsign}) this explains qualitatively why the lower eigenvalues converge with alternate signs, whereas the higher ones display an almost monotonic behaviour.  

\begin{center}
\begin{figure}[!t]
\begin{center}
\includegraphics[width=0.7\textwidth]{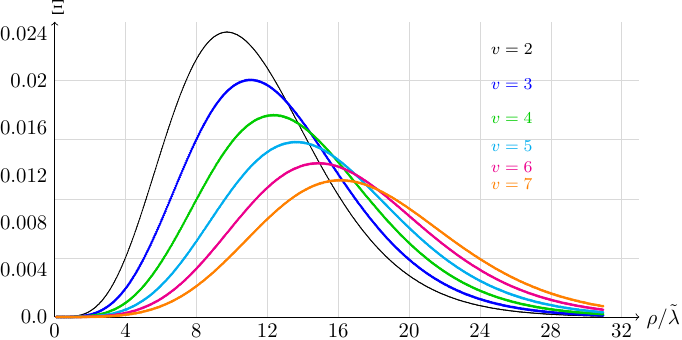}
\vskip 0.2cm
\caption{\small Positivity of the linear combination $\Xi(\rho/\tilde\lambda)$ of
  $\chi^2$--ratios, eq.~(\ref{eq:xieq}).\label{fig:Fineq}}       
\end{center}
\vskip 0.0cm
\end{figure}
\end{center}

\section{Perturbative estimators vs. the iterative one}

In the last part of the paper we introduce some statistical noise. So far we have studied the eigenvalue reconstruction under the hypothesis that $\mu = \tau_\rho\cdot\lambda$ represents the exact truncated counterpart of some $\lambda\in\RR^v$. This is rather unusual in most applications. In general, $\mu$ is not the result of an exact truncation. From now on we assume it to be estimated from a representative sample $\cP_N = \{x^{(k)}\}_{k=1}^N$ of $X\sim\cN(0,\Sigma)$ with finite size $N$. As usual in sample space, we regard the observations $x^{(k)}$ as realizations of i.i.d. stochastic variables $X^{(k)}\sim\cN(0,\Sigma)$, $k=1,\ldots,N$. Sample estimates of $\mu$ are performed as follows. A certain subset of $M<N$ elements of $\cP_N$ falls within $\cB_v(\rho)$, with the fraction $M/N$ fulfilling $\lim_{N\to\infty} M/N = \alpha(\rho;\lambda)$. From this subset we measure ${\frak S}_\cB$  via the estimator  
\begin{align}
(\hat{\frak S}_\cB)_{ij} & = \frac{1}{M-1}\sum_{k=1}^N (x^{(k)}-\hat x)_i\cdot
(x^{(k)}-\hat x)_j \cdot \I_{\cB_v(\rho)}(x^{(k)})\,,\\[0.0ex]
\hat x_i & = \frac{1}{M}\sum_{k=1}^N  x^{(k)}_i\cdot \I_{\cB_v(\rho)}(x^{(k)})\,,
\end{align}
with $\I_{\cB_v(\rho)}(\,\cdot\,)$ denoting the characteristic function of $\cB_v(\rho)$. We showed in ref.~\cite{palombi4} that $\hat x_i$ is an unbiased estimator of the truncated mean, while $\hat{\frak S}_\cB$ is asymptotically unbiased. The vector $\hat\mu$ of the eigenvalues of $\hat{\frak S}_\cB$ represents our sample estimate of $\mu$. We use $\hat\mu$ as an input parameter for the iterative reconstruction of $\lambda$, conditioned to $\hat\mu\in\cD(\tau_\rho^{-1})$, and for the perturbative reconstruction of $\lambda$, conditioned to $\hat\mu\in\cD(\cT_\rho^{-1})$. We refer the reader to ref.~\cite{palombi4} for a discussion of the failure probabilities
\begin{align}
  p_\text{fail}(\rho,\Sigma,N) & = \mathds{P}\left[\hat\mu\notin\cD(\tau_\rho^{-1})\,|\,X^{(k)}\sim\cN(0,\Sigma)\,,\ k=1,\ldots, N \right]\,,\\[0.0ex]
  q_\text{fail}(\rho,\Sigma,N) & = \mathds{P}\left[\hat\mu\notin\cD(\cT_\rho^{-1})\,|\,X^{(k)}\sim\cN(0,\Sigma)\,,\ k=1,\ldots, N \right].
\end{align}
The variable $\hat\mu$ can be interpreted in its turn as the realization of a stochastic variable in sample space. It thus makes sense to pose the question of what statistical properties characterize the stochastic variable $\hat\lambda = \tau_\rho^{-1}\cdot\hat\mu$, conditioned to $\hat\mu\in\cD(\tau_\rho^{-1})$, and its perturbative approximations $\hat\lambdabar$\raisebox{-0.1ex}{$^{(k)}$}~($k=1,\ldots,4$), conditioned to $\hat\mu\in\cD(\cT_\rho^{-1})$.

Finding analytic relations between the expectation in sample space of polynomial functions of $\hat\mu$ and analogous functions of $\hat\lambda$ is non-trivial, since no analytic representation of $\tau_\rho^{-1}$ is given. This task goes beyond the aims of the present paper. Here, we adopt a pragmatic approach where we limit ourselves to simulations with a specific choice of $\Sigma$. In particular, we assume that $\cP_N$ distributes normally with $\Sigma=\diag(\lambdaEX)$ and $\lambdaEX$ as introduced in sect.~4.  In our study, we choose $N=200,\ 250,\ \ldots,\ 2000$; for each value of $N$, we generate about 5000 normal populations; for each of them, we consider Euclidean balls with $\rho=4.0,\ 6.0,\ldots,\ 40.0$ and for each pair $(\rho,N)$ we measure bias and variance of $\hat\lambda$ and $\hat\lambdabar$\raisebox{-0.1ex}{$^{(k)}$}. In Fig.~\ref{fig:biasvar} we report the results we obtained for $\rho=6.0$ (corresponding to weak truncation with $\alpha(6.0;\lambdaEX)\simeq 0.844$)\footnote{Statistical errors of the sample estimate of the variances have been computed according to the general formula for the standard error ${\rm se}(\var(\hat\lambda_k)) = \widehat{\var(\hat\lambda_k)}\cdot[2/(N-1)+\hat\kappa/N]$, where $\widehat{\var(\hat\lambda_k)}$ is the sample estimate of $\var(\hat\lambda_k)$ and  $\hat\kappa$ is the sample excess kurtosis of the distribution of $\var(\hat\lambda_k)$.\phantom{$\widehat{\var(\hat\lambda_k)}$}}. From the plots on the left we notice that
\begin{itemize}[itemsep=0.0em]
\item[{\it i})]{the bias of $\hat\lambdabar$\raisebox{-0.1ex}{$^{(k)}$} is weakly sensitive to $N$  for all $k$'s; it converges asymptotically to an intrinsic perturbative bias with finite size corrections proportional to~$1/N$;}
\item[{\it ii})]{$\hat\lambda$ is slightly biased at finite $N$ and asymptotically unbiased; convergence to the asymptotic limit is again reached linearly in $1/N$.}  
\end{itemize}
Similarly, from the plots on the right we observe that
\begin{itemize}[itemsep=0.0em]
\item[{\it iii})]{all variances vanish asymptotically and have finite size corrections proportional to $1/N$;}
\item[{\it iv})]{the variance of the higher eigenvalues ($\lambdabar^{(k)}_3,\lambdabar^{(k)}_4$) increases at fixed $N$ as we increase the order of truncation of the perturbative series;}
\item[{\it v})]{the variance of the lower eigenvalues ($\lambdabar^{(k)}_1,\lambdabar^{(k)}_2$) decreases at fixed $N$ as we increase the order of truncation of the perturbative series;}
\item[{\it vi})]{the iterative estimator of the higher eigenvalues ($\hat\lambda_3,\hat\lambda_4$) has a higher variance than all its perturbative approximations, while the iterative estimator of the lower eigenvalues ($\hat\lambda_1,\hat\lambda_2$) has a lower variance than all its perturbative approximations.}
\end{itemize}

The variance plots illustrate the potential usefulness of the perturbative estimators. The reconstruction of the higher eigenvalues achieved from the iterative algorithm is noisy for moderately small $N$ (say -- $N\lesssim 400$). Perturbative estimators allow to control the variance. The price to pay for this is the introduction of a non-vanishing asymptotic bias. Depending on the specific context, there may be an optimal choice for the order of truncation of the perturbative series, which guarantees acceptable values of both bias and variance.

We find qualitatively similar results for other values of $\rho$: both the asymptotic biases and the slopes of the variances decrease as $\rho$ increases, as intuitively expected.

An important result of our simulations is inferred upon relating the variance of the reconstructed eigenvalues to that of the truncated ones. We observe that since $\lambda = \tau_\rho^{-1}\cdot\mu$ is a vector relation, each of the reconstructed eigenvalues $\lambda_i=\lambda_i(\mu_1,\ldots,\mu_v)$ depends upon all the components of $\mu$. It follows that also ${\rm var}(\lambda_i)$ is a function of all the components of $\mu$. Nevertheless, $\lambda_i$ depends weakly on $\mu_k$ for $k\ne i$. Therefore, it makes sense to examine how $\var(\hat\lambda_i)$ relates to $\var(\hat\mu_i)$. An example of such dependence is shown in Fig.~\ref{fig:varvar} for $i=1,4$, corresponding respectively to the lowest and highest components of $\lambdaEX$. We note that the variances are linearly related, except for weak quadratic corrections observed for $\var(\hat\mu_4)\simeq 1.0\times 10^{-2}$. We also observe that the variance of the iterative estimator of the lowest eigenvalue is minimal and that of the highest one is maximal. The slopes observed for the highest eigenvalue $\hat\lambda_4$ are remarkable. By comparing the scales of the $x$- and $y$-axis we recognize that a huge inflation of the variance occurs as a result of applying $\tau_\rho^{-1}$ to $\hat\mu$. Numerical simulations signal the existence of such amplification phenomena, which are ultimately due to the unboundedness of $\tau_\rho^{-1}$. In practical situations the exact reconstruction of the highest eigenvalue may be critical. In such cases, the adoption of perturbative estimators in place of the iterative one may represent a viable solution.

\section{Conclusions}

In this paper we have explored a perturbative approach to the reconstruction of a normal covariance matrix $\Sigma$ from a spherically truncated counterpart ${\frak S}_\cB$. Since $\Sigma$ and ${\frak S}_\cB$ are simultaneously diagonalized, the reconstruction problem concerns only their eigenvalues. After collecting all the ingredients needed for the algebraic implementation of the perturbative expansion of the reconstruction operator, we have examined the general structure of the perturbative series and some practical aspects related to the calculation of the first few perturbative coefficients. We provide formulae for the reconstruction of the eigenvalues of $\Sigma$ up to the fourth perturbative order as well as \Maple programs to further improve the approximation. 

From a theoretical standpoint, the perturbative method is meant to complement the fixed--point iterative algorithm proposed by us in ref.~\cite{palombi4} in cases where the covariance reconstruction is ill-defined or the iterative algorithm is inefficient. Such cases occur when either ${\frak S}_\cB$ is affected by stochastic noise, the squared truncation radius $\rho$ is comparable or less than the lowest eigenvalue of $\Sigma$, or the number of dimensions $v$ is large. The ill-posedness of the reconstruction problem emerges when, due to statistical fluctuations, the eigenvalues of ${\frak S}_\cB$ lie outside the domain of the reconstruction operator. In this case the perturbative approach provides a regularization with respect to the existence of a solution. Instead, in cases of small $\rho$ or large $v$, the inefficiency of the iterative algorithm consists in slow convergence speed.

Another weakness of the iterative algorithm emerges when the eigenvalue reconstruction is performed from  statistically poor samples of the eigenvalues of ${\frak S}_\cB$. We have shown that the statistical noise of these is inflated by the application of the reconstruction operator, thus producing large fluctuations of the higher components of the reconstructed eigenvalues. Perturbation theory offers the possibility to control the variance and stabilize the reconstruction by properly choosing the order of truncation of the perturbative series. The price to pay when replacing the iterative estimator with perturbative approximations is the introduction of an asymptotic non-vanishing bias. It is possible to adopt mixed approaches, where the lower eigenvalues are reconstructed via the iterative algorithm while the higher ones are obtained from perturbation theory. 

\section*{Acknowledgements}

The computing resources used for our numerical study and the related technical support have been partly provided by the CRESCO/ENEAGRID High Performance Computing infrastructure and its staff \cite{Ponti}. CRESCO ({\color{red} C}omputational {\color{red} RES}earch centre on {\color{red} CO}mplex systems) is funded by ENEA and by Italian and European research programmes.

\begin{center}
\begin{landscape}
\begin{figure}[!t]
\begin{center}
\includegraphics[width=1.15\textwidth]{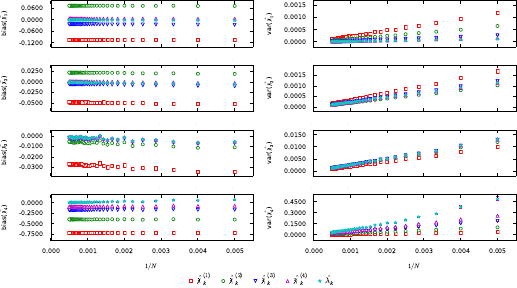}
\vskip 0.4cm
\caption{\small Bias and variance of the iterative and perturbative estimators vs. the inverse of the population size for $\rho=6.0$. Normal populations have been generated with $\Sigma=\diag(\lambdaEX)$.\label{fig:biasvar}}     
\end{center}
\end{figure}
\end{landscape}
\end{center}
\begin{center}
\begin{figure}[!t]
\begin{center}
\includegraphics[width=0.7\textwidth]{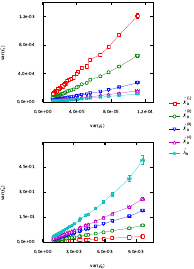}
\vskip 0.4cm
\caption{\small Variances of the iterative and perturbative estimators of the lowest and highest reconstructed eigenvalue vs. the variance of the corresponding truncated eigenvalue. Normal populations have been generated with $\Sigma=\diag(\lambdaEX)$.\label{fig:varvar}}    
\end{center}
\end{figure}
\end{center}

\begin{appendices}

  \clearpage

\section{Injectivity of the operator $\tau_\rho$}

We let $\lambda',\lambda''\in \dR^v_+$ denote two variance vectors. For $\rho\in\dR_+$, we also let $\mu'_k = \lambda'_k\frac{\alpha_k}{\alpha}(\rho;\lambda')$ and $\mu_k''=\lambda_k''\frac{\alpha_k}{\alpha}(\rho;\lambda'')$ for $k=1,\ldots,v$. We want to show that if $\mu'=\mu''$, then $\lambda' = \lambda''$. In consideration of the Fundamental Theorem of Calculus for line integrals, under the assumption that $\mu'=\mu''$, we have
\begin{equation}
0 = \mu'_k - \mu_k'' = \lambda'_k\frac{\alpha_k}{\alpha}(\rho;\lambda') -
\lambda_k''\frac{\alpha_k}{\alpha}(\rho;\lambda'') = \sum_{\ell=1}^v \left[\int_0^1\rd t\
J_{k\ell}\left(\rho;\lambda'' + t\left(\lambda'-\lambda''\right)\right)\right](\lambda'_\ell - \lambda''_\ell)\,,
\label{eq:fixedpointid}
\end{equation}
where $J$ denotes the Jacobian of $\tau_\rho$, having matrix elements
\begin{equation}
J_{k\ell}(\rho;\lambda)
= \partial_k\left(\lambda_\ell\frac{\alpha_\ell}{\alpha}(\rho;\lambda)\right) =
\frac{1}{2}\frac{\lambda_\ell}{\lambda_k}\left(\frac{\alpha_{k\ell}}{\alpha} - \frac{\alpha_k\alpha_\ell}{\alpha^2}\right) = \left[\Lambda^{-1}\,\Omega(\rho;\lambda)\,\Lambda\right]_{k\ell}\,,
\end{equation}
with $\Omega_{k\ell}\equiv (1/2)(\alpha_{k\ell}/\alpha - \alpha_k\alpha_\ell/\alpha^2)$. From eq.~(23) of ref.~\cite{palombi4} we know that $\Omega(\rho;\lambda) = \{\Omega_{k\ell}(\rho;\lambda)\}_{k,\ell=1}^v$ is the covariance matrix of the square components of $X\sim\cN_v(0,\text{diag}(\lambda))$ under spherical truncation with square radius $\rho$. As such, $\Omega(\rho;\lambda)$ is symmetric and positive definite. Explicitly, we have
\begin{align}
\Omega_{k\ell} & =\,
\frac{1}{2\lambda_k\lambda_\ell}\cov\left(X_k^2,X_\ell^2\,|\,X\in\cB_v(\rho)\right)
\nonumber\\[1.0ex] 
& = \, \frac{1}{2\lambda_k\lambda_\ell}\E\biggl[\biggl(X_k^2-\E\left[X_k^2\,|\,X\in\cB_v(\rho)\right]\biggr)\biggl(X_\ell^2-\E[X_\ell^2\,|\,X\in\cB_v(\rho)]\biggr)\,\biggr|\,X\in\cB_v(\rho)\biggr]\,.
\end{align}
On setting $Z_k = (X_k^2-\E[X_k^2\,|\,X\in\cB_v(\rho)])/\sqrt{2}\lambda_k$, we represent
$\Omega$ as $\Omega = \E[Z\trans{Z}\,|\,X\in\cB_v(\rho)]$. If $x\in\RR^v$ is not the null vector,
then $\trans{x}\Omega x = \E[\trans{x}Z\trans{Z}x\,|\,X\in\cB_v(\rho)] =
\E[(\trans{x}Z)^2\,|\,X\in\cB_v(\rho)]>0$. Moreover, the eigenvalues of $\Omega$ 
fulfill the secular equation 
\begin{equation}
0=\det(\Omega - \phi\mathds{I}_v) = \det[\Lambda^{-1}(\Omega -
\phi\mathds{I}_v)\Lambda] = \det(\Lambda^{-1}\Omega\Lambda-\phi\mathds{I}_v) = \det(J-\phi\mathds{I}_v)\,.
\end{equation}
It follows that $J$ is positive definite (though it is not
symmetric). Since the sum of positive definite matrices is positive definite,
we conclude that $\int_0^1\rd t\ J\left(\rho;\lambda'' +
  t\left(\lambda'-\lambda''\right)\right)$ is positive definite too. As such, it
is non--singular. Hence, we conclude from eq.~(\ref{eq:fixedpointid}) that
$\lambda'=\lambda''$.  

\section{Domain of the operator $\tau_\rho^{-1}$}

In ref.~\cite{palombi4}, we proved the following two properties of the truncation operator:

\begin{prop}[monotonicities]
\label{prop:mon}
Let $\lambda\in\dR_+^v$ denote a variance vector and $\lambda_{(k)}\equiv\{\lambda_i\}_{i=1,\ldots,v}^{i\ne k}$ the set of variances without $\lambda_k$. For $\rho\in\dR_+$, the variances truncated at $\rho$ fulfill the following properties:  
\begin{itemize}
\item[{\rm (}$p_1${\rm )}]{\ $\lambda_k\dfrac{\alpha_k}{\alpha}(\rho;\lambda)$ is a
  monotonic increasing function of $\lambda_k$ for fixed $\rho$ and $\lambda_{(k)}$;}
\item[{\rm (}$p_2${\rm )}]{\ $\lambda_k\dfrac{\alpha_k}{\alpha}(\rho;\lambda)$ is a monotonic
  decreasing function of $\lambda_i$ $(i\ne k)$ for fixed $\rho$ and $\lambda_{(i)}$,}
\end{itemize}
\end{prop}

\begin{prop}[variance ordering]
  \label{prop:VarOrd}
Let $\lambda\in\dR_+^v$ denote a variance vector and, for $\rho\in \dR_+$, let $\mu\in\dR_+^v$ be the vector of variances truncated at $\rho$. If $\lambda_{1}\le \lambda_{2}\le\ldots\le\lambda_{v}$, then $\mu_{1}\le \mu_{2}\le\ldots\le\mu_{v}$. 
\end{prop}
\noindent We shall not repeat the proofs here. Prop.~\ref{prop:VarOrd} allows us to split $\text{Im}(\tau_\rho)$ into non-overlapping sectors. Specifically, we let
\begin{align}
  \Im(\tau_\rho) & = \biggl\{\mu\in\dR^v_+:\ \mu_k = \lambda_k \frac{\alpha_k}{\alpha}(\rho;\lambda) \ \text{ for } \  k = 1,\ldots,v \ \nonumber\\[-2.0ex]
  & \hskip 2.4cm \text{ and for some } \lambda\in\dR_+^v \text{ with } \lambda_1\le\ldots\le\lambda_v\biggr\}\,.
\end{align}
Accordingly, we have
\begin{equation}
  \text{Im}(\tau_\rho) = \bigcup_{\sigma\in S_v}\sigma\circ \Im(\tau_\rho) = \{\mu: \mu = \sigma\circ \mu_0 \text{ for } \mu_0\in\Im(\tau_\rho) \text{ and } \sigma\in S_v \}\,,
\end{equation}
with $S_v$ being the set of permutations of $v$ elements. Hence, we can focus on $\Im(\tau_\rho)$. To characterize its boundary $\partial\Im(\tau_\rho)$, we use Prop.~\ref{prop:mon}. Specifically, we look for the limit values of $\mu$ as $\lambda\to\infty$ along a sequence of properly chosen directions. To this aim, we must respect the increasing order of the components of $\lambda$ for $\mu\in\Im(\tau_\rho)$. For instance, we cannot let $\lambda_{v-1}\to\infty$ while keeping $\lambda_v$ fixed. To overcome the problem, we introduce the \emph{added truncated moments}
\begin{equation}
  \nu_k(\rho;\lambda_1,\ldots,\lambda_k) = \lambda_k\sum_{i=k}^v \frac{\alpha_i}{\alpha}(\rho;\lambda)\biggr|_{\lambda_{k+1}=\ldots=\lambda_v=\lambda_k}\,,\qquad k=1,\ldots,v\,.
\end{equation}
As a consequence of Prop.~\ref{prop:mon}, we can show that
\begin{itemize}
\item[{\rm (}$p_3${\rm )}]{\ $\nu_k(\rho;\lambda_1,\ldots,\lambda_k)$ is a
  monotonic increasing function of $\lambda_k$ for fixed $\rho$ and $\lambda_i$ ($i\ne k$);}
\item[{\rm (}$p_4${\rm )}]{\ $\nu_k(\rho;\lambda_1,\ldots,\lambda_k)$ is a monotonic
  decreasing function of $\lambda_i$ $(i<k)$ for fixed $\rho$ and $\lambda_j$ ($j\ne i$),}
\end{itemize}
Indeed, we have
\begin{equation}
  \nu_k(\rho;\lambda_1,\ldots,\lambda_k) = \frac{\int_{\cB_v(\rho)}\rd^vx\,\left(\sum_{i=k}^vx_i^2\right)\prod_{i=1}^{k-1}\delta(x_i,\lambda_i)\prod_{i=k}^{v}\delta(x_i,\lambda_k)}{\int_{\cB_v(\rho)}\rd^vx\prod_{i=1}^{k-1}\delta(x_i,\lambda_i)\prod_{i=k}^{v}\delta(x_i,\lambda_k)}\,.
\end{equation}
Differentiating $\nu_k$ under the integral sign yields
\begin{equation}
  \partial_k\nu_k(\rho;\lambda_1,\ldots,\lambda_k) = \frac{1}{2\lambda_k^2}\text{var}[X_k^2+\ldots+X_v^2\,|\,X\in\cB_v(\rho)]\ge 0\,.
\end{equation}
On the other hand, for $i<k$ we have
\begin{equation}
  \partial_i\nu_k(\rho;\lambda_1,\ldots,\lambda_k) = \sum_{j=k}^v\partial_i\left[\lambda_j\frac{\alpha_j}{\alpha}(\rho;\lambda)\biggr|_{\lambda_{k+1}=\ldots=\lambda_v=\lambda_k}\right]\le 0\,,
\end{equation}
because each term of the sum is negative in view of ($p_2$). It is also important to notice that all terms in $\nu_k$ are equal by symmetry, hence
\begin{equation}
  \nu_k(\rho;\lambda_1,\ldots,\lambda_k) = (v-k+1)\mu_j(\rho;\lambda)|_{\lambda_{k+1}=\ldots=\lambda_v=\lambda_{k}}\,,\qquad j=k,k+1,\ldots,v\,,
\end{equation}
or, equivalently, 
\begin{equation}
  \mu_j(\rho;\lambda)|_{\lambda_{k+1}=\ldots=\lambda_v=\lambda_{k}} = \frac{1}{v-k+1}\nu_k(\rho;\lambda_1,\ldots,\lambda_k)\,,\qquad j=k,k+1,\ldots,v\,.
  \label{eq:addedone}
\end{equation}
Moreover, we can calculate exactly the limit of $\nu_k$ as $\lambda_k\to\infty$ and $\lambda_1,\ldots,\lambda_{k-1}\to 0$. Using spherical coordinates, we find
\begin{align}
  \tilde \nu_k(\rho) & = \lim_{\lambda_1,\ldots,\lambda_{k-1}\to 0}\lim_{\lambda_k\to\infty}\nu_k(\rho;\lambda_1,\ldots,\lambda_k) = \frac{\int_{\sum_{i=k}^vx_i^2<\rho}\rd x_k\ldots\rd x_v \sum_{i=k}^{v}x_i^2}{\int_{\sum_{i=k}^vx_i^2<\rho}\rd x_k\ldots\rd x_v}\nonumber\\[2.0ex]
  & = \frac{\int_0^{\sqrt{\rho}}\rd r\,r^{v-k+2}}{\int_0^{\sqrt{\rho}}\rd r\, r^{v-k}} = \frac{v-k+1}{v-k+3}\,\rho\,.
  \label{eq:addedtwo}
\end{align}
In particular, we have
\begin{equation}
  \tilde \nu_v(\rho) = \frac{1}{3}\rho\,, \quad \tilde \nu_{v-1}(\rho) = \frac{2}{4}\rho\,, \quad \ldots \quad \,, \quad  \tilde \nu_1(\rho) = \frac{v}{v+2}\,\rho\,.
  \label{eq:addedthree}
\end{equation}
From Eqs.~(\ref{eq:addedone})-(\ref{eq:addedthree}), we conclude that all points
\begin{equation}
  \arraycolsep=1.4pt\def\arraystretch{1.4}
  \begin{array}{lccccccccc}
  P_v     & = & \bigl( & 0, & \ldots, & 0, & 0, & 0, & \rho/3 & \bigr)\,,\\
  P_{v-1} & = & \bigl( & 0, & \ldots, & 0, & 0, & \rho/4, & \rho/4 & \bigr)\,,\\
  P_{v-2} & = & \bigl( & 0, & \ldots, & 0, & \rho/5, & \rho/5, & \rho/5 & \bigr)\,,\\[1ex]
  & \vdots & &&&&&&&\\[0ex]
  P_{1} & = & \bigl( & \rho/(v+2), & \ldots, & \rho/(v+2), & \rho/(v+2), & \rho/(v+2), & \rho/(v+2) & \bigr)\,, \\
  \end{array}
  \label{eq:cusps}
\end{equation}
belong to $\partial\Im(\tau_\rho)$, i.e. they fulfill $\sum_{\ell=1}^{v-1}(P_k)_\ell + 3(P_k)_v = \rho$, $k=1,\ldots,v$. In particular, for $v=3$ points $P_1 = \{\rho/5,\rho/5,\rho/5\}$, $P_2 = \{0,\rho/4,\rho/4\}$, $P_3 = \{0,0,\rho/3\}$ and those obtained by permuting their components in all possible ways yield the cusps of $\cD(\tau_\rho^{-1})$ in Fig.~\ref{fig:defdoms}(left).  This observation suggests that $\partial\Im(\tau_\rho)$ can be also written as the convex hull of Eq.~(\ref{eq:cusps}), namely
\begin{equation}
  \partial\Im(\tau_\rho) = \text{Conv}\left(\{P_1,\ldots,P_v\}\right) = \left\{\sum_{i=1}^v a_i P_i\, \bigl|\,(\forall i:a_i\ge 0) \text{ and } \sum_{i=1}^va_i = 1\right\}\,.
  \label{eq:convhull}
\end{equation}
Unfortunately, we lack a formal proof of Eq.~(\ref{eq:convhull}). The argument we used to calculate limit values for $\mu_k$ fails as soon as we let $\lambda\to\infty$ along any other direction than $\lambda_1,\ldots,\lambda_{k-1}\to 0$, $\lambda_k=\lambda_{k+1}=\ldots=\lambda_v\to\infty$, for $k=1,\ldots,v$, due to symmetry breaking.

\section{\Maple code}

\subsection{Code block 1: the coefficient $\Delta_{i_1\ldots i_n}$}

\vskip 0.3cm
{\footnotesize
\begin{verbatim}
  # Delta coefficient
  # -----------------

  Delta := proc()
    local SortArgs,V,ActCtr,NxtCtr,Res,k:
    for k from 1 to _npassed do
      if not type(_passed[k],'nonnegint') then
        return 'procname(_passed)':
      end if:
    end do:
    SortArgs := sort([_passed[1.._npassed]]):
    V := Vector(_npassed):
    V[1] := 1:  
    ActCtr := 1:
    NxtCtr := 2:
    for k from 1 to (_npassed-1) do
      if SortArgs[NxtCtr] = SortArgs[ActCtr] then
        V[ActCtr] := V[ActCtr]+1:
      else
        ActCtr := NxtCtr:
        V[ActCtr] := 1:
      end if:
      NxtCtr := NxtCtr+1:
    end do:
    Res := 1:
    for k from 1 to _npassed do
      Res := Res*(doublefactorial(2*V[k]-1)):
    end do:
    return Res:
  end:
\end{verbatim}
}

\subsection{Code block 2: perturbative expansion of eq.~(\ref{eq:truncspectrum})}

\vskip 0.3cm
{\footnotesize
\begin{verbatim}
  # Nested sequence
  # ---------------
  
  NestSeq := proc(TheEq,v::nonnegint,niter::nonnegint)
    if niter = 0 then
      eval(TheEq):
    else
      seq(eval(NestSeq(TheEq,v,niter-1)),
          cat('r', niter) = 1..v):
    end if
  end proc:
  
  # Derivatives of Gaussian Integrals
  # ---------------------------------
  
  DerAlpha := proc()
    global v,Delta:
    local m,Fact1,Fact2:
    m := _npassed:
    Fact1 := Delta(_passed[1.._npassed])/(2*l[0])^m:
    Fact2 := add((-1)^(m-j)*binomial(m,j)*F[v+2*j],j=0..m):
    return Fact1*Fact2:
  end proc:
  
  DerAlphak := proc()
    global v,Delta:
    local m,Fact1,Fact2:
    m := _npassed-1:
    Fact1 := Delta(_passed[1.._npassed])/(2*l[0])^m:
    Fact2 := add((-1)^(m-j)*binomial(m,j)*F[v+2*(j+1)],j=0..m):
    return Fact1*Fact2:
  end proc:
  
  # Function arguments
  # ------------------
  
  lam := Vector(v):
  for k from 1 to v do
    lam[k] := add(l[j,k]*epsilon^j,j=0..n):
  end do:
  lam := seq(lam[k],k=1..v):
  lam0 := seq(l[0,k],k=1..v):
  
  # Integral ratio
  # --------------
  
  R := proc(j)
    global lam:
    return alpha[j](lam)/alpha(lam):
  end:
  
  # Taylor expansion of the map
  # ---------------------------
  
  for j from 1 to n do
    for k from 1 to v do
      Rk := convert(taylor(R(k),epsilon=0,n+1),polynom):
      h[j,k] := expand(coeff(lam[k]*Rk,epsilon,j)):
    end do:
  end do:
  
  # Evaluation conditions
  # ---------------------
  
  C0A := alpha(lam0)=F[v]:
  C0Ak := seq(alpha[k](lam0)=F[v+2],k=1..v):
  
  C1A := NestSeq(D[r1](alpha)(lam0)=DerAlpha(r1),v,1):
  C1Ak := NestSeq(D[r1](alpha[r2])(lam0)=DerAlphak(r1,r2),v,2):
  
  C2A := NestSeq(D[r1,r2](alpha)(lam0)=DerAlpha(r1,r2),v,2):
  C2Ak := NestSeq(D[r1,r2](alpha[r3])(lam0)=DerAlphak(r1,r2,r3),v,3):
  
  C3A := NestSeq(D[r1,r2,r3](alpha)(lam0)=DerAlpha(r1,r2,r3),v,3):
  C3Ak := NestSeq(D[r1,r2,r3](alpha[r4])(lam0)=DerAlphak(r1,r2,r3,r4),v,4):
  
  C4A := NestSeq(D[r1,r2,r3,r4](alpha)(lam0)=DerAlpha(r1,r2,r3,r4),v,4):
  C4Ak := NestSeq(D[r1,r2,r3,r4](alpha[r5])(lam0)=DerAlphak(r1,r2,r3,r4,r5),v,5):
  
  CArg0 := seq(l[0,k]=l[0],k=1..v):
  
  # Evaluations
  # -----------
  
  for j from 1 to n do
    for k from 1 to v do
      h00[j,k] := expand(eval(h[j,k],[C1A])):
      h01[j,k] := expand(eval(h00[j,k],[C2A])):
      h02[j,k] := expand(eval(h01[j,k],[C3A])):
      h03[j,k] := expand(eval(h02[j,k],[C4A])):
      h04[j,k] := expand(eval(h03[j,k],[C1Ak])):
      h05[j,k] := expand(eval(h04[j,k],[C2Ak])):
      h06[j,k] := expand(eval(h05[j,k],[C3Ak])):
      h07[j,k] := expand(eval(h06[j,k],[C4Ak])):
      h08[j,k] := expand(eval(h07[j,k],[C0A])):
      h09[j,k] := expand(eval(h08[j,k],[C0Ak])):
      h10[j,k] := expand(eval(h09[j,k],[CArg0])):
    end do:
  end do:
\end{verbatim}
}

\subsection{Code block 3: extraction of $\gamma^{(n,m)}_{k_1\ldots k_{n+1}}$}

\vskip 0.3cm
{\footnotesize
\begin{verbatim}
with(LinearAlgebra):
with(RandomTools):

v := 6:

# O-structure matrix
# ------------------

zeta1 := add(l[1,k],k=1..v):
zeta2 := add(l[2,k],k=1..v):
zeta11 := add(l[1,k]^2,k=1..v):
zeta12 := add(l[1,k]*l[2,k],k=1..v):
zeta111 := add(l[1,k]^3,k=1..v):

S3matrix := Matrix(v,v):

for k from 1 to v do
  S3matrix[k,1] := l[1,k]^3:
  S3matrix[k,2] := expand(l[1,k]*zeta11):
  S3matrix[k,3] := expand(zeta111):
  S3matrix[k,4] := l[0]*l[1,k]*l[2,k]:
  S3matrix[k,5] := expand(l[0]*l[1,k]*zeta2):
  S3matrix[k,6] := expand(l[0]*zeta12):
end do:

# Jacobian matrix
# ---------------

Jmatrix := Matrix(v,v):

for k1 from 1 to v do
  for k2 from 1 to v do
    Jmatrix[k1,k2] := (1/2)*(Delta(k1,k2)*F[v+4]/F[v] - F[v+2]^2/F[v]^2):
  end do:
end do:

# Terms to be removed by hand
# ---------------------------

V := Vector(v):

for j from 1 to v do
  V[j] := 0:
  for k from 1 to v do:
    V[j] := V[j] + Jmatrix[j,k]*l[3,k]:
  end do:
end do:

# Randomized Linear system
# ------------------------

C := Vector(v):

for j from 1 to v do
  lincond0 := l[0]=Generate(float(range=0..1,'method=uniform')):
  linvals1 := seq(Generate(float(range=0..1,'method=uniform')),m=1..v-1):
  lincond1 := seq(l[1,m]=linvals1[m],m=1..v-1):
  lincond1 := lincond1,l[1,v]=-add(linvals1[m],m=1..v-1):
  lincond2 := seq(l[2,m]=Generate(float(range = 0..1,'method=uniform')),m=1..v):
  for m from 1 to v do
    S3matrix[j,m] := eval(S3matrix[j,m],[lincond0,lincond1,lincond2]):
  end do:

  r := expand(h10[3,j] - V[j]):
  s := expand(eval(coeff(r,F[v+6]),F[v+2]=0)):
  C[j] := eval((l[0]^2)*F[v]*s,[lincond0,lincond1,lincond2]):
end do:

# (-1) x Gamma coefficients
# -------------------------

Gcoefs := LinearSolve(S3matrix,C):
\end{verbatim}
}

\end{appendices}

\bibliographystyle{hunsrt}
\bibliography{main}

\end{document}